\documentclass[12pt,a4paper]{article}
\usepackage{amssymb}

\usepackage{amsmath}

\newtheorem{theorem}{Theorem}

\newtheorem{corollary}[theorem]{Corollary}

\newtheorem{definition/proposition}[theorem]{Definition/Proposition}
\newtheorem{definition}[theorem]{Definition}

\newtheorem{lemma}[theorem]{Lemma}

\newtheorem{proposition}[theorem]{Proposition}

\input{tcilatex}

\begin{document}

\date{}
\title{Fragmentation processes with an initial mass converging to infinity}
\author{B\'{e}n\'{e}dicte Haas\thanks{%
Department of Statistics, University of Oxford, 1 South Parks
Road, Oxford OX1 3TG, UK. E-mail: haas@ceremade.dauphine.fr.
Research supported in part by EPSRC GR/T26368.} } \maketitle

\begin{abstract}
We consider a family of fragmentation processes where the rate at which a
particle splits is proportional to a function of its mass. Let $%
F_{1}^{(m)}(t),F_{2}^{(m)}(t),...$ denote the decreasing rearrangement of
the masses present at time $t$ in a such process, starting from an initial
mass $m$. Let then $m\rightarrow \infty $. Under an assumption of regular
variation type on the dynamics of the fragmentation, we prove that the
sequence $(F_{2}^{(m)},F_{3}^{(m)},...)$ converges in distribution, with
respect to the Skorohod topology, to a fragmentation with immigration
process. This holds jointly with the convergence of $m-F_{1}^{(m)}$ to a
stable subordinator. A continuum random tree counterpart of this result is
also given: the continuum random tree describing the genealogy of a
self-similar fragmentation satisfying the required assumption and starting
from a mass converging to $\infty $ will converge to a tree with a spine
coding a fragmentation with immigration.
\end{abstract}

\ \


\textbf{A.M.S. Classification. }60J25, 60F05.

\section{Introduction and main results}

We consider Markovian models for the evolution of systems of particles that
undergo splitting, so that each particle evolves independently of others
with a splitting rate proportional to a function of its mass. In \cite%
{bertsfrag02}, Bertoin obtains such fragmentation model with some
self-similarity property by cutting the Brownian Continuum Random Tree (CRT)
of Aldous \cite{AldousCRT1},\cite{Aldous CRT3} as follows: for all $t\geq 0$%
, remove all the vertices of the Brownian CRT that are located under height $%
t$ and consider the connected components of the remaining vertices. Next,
set $F^{B_{r},(1)}(t):=(F_{1}^{B_{r},(1)}(t),F_{2}^{B_{r},(1)}(t),...)$ for
the decreasing sequence of masses of these connected components: $%
F^{B_{r},(1)}$ is then a fragmentation process starting from $(1,0,...)$
where fragments split with a rate proportional to their mass to the power $%
-1/2.$

On the other hand, Aldous \cite{AldousCRT1} shows that the Brownian CRT
rescaled by a factor $1/\varepsilon $ converges in distribution to an
infinite CRT composed by an infinite baseline $\left[ 0,\infty \right) $ on
which are attached compact CRT's distributed, up to a scaling factor, as the
Brownian CRT. In terms of fragmentations, his result implies that%
\begin{equation*}
\varepsilon ^{-2}(F_{2}^{B_{r},(1)}(\varepsilon \cdot
),F_{3}^{B_{r},(1)}(\varepsilon \cdot ),...)\overset{\text{law}}{\rightarrow
}FI^{B_{r}}\text{ as }\varepsilon \rightarrow 0
\end{equation*}%
where $FI^{B_{r}}$ is some fragmentation with immigration process
constructed from the infinite Brownian CRT of Aldous. Equivalently, if $%
F^{B_{r},(m)}$ denotes the Brownian fragmentation starting from $(m,0,...),$%
\begin{equation*}
(F_{2}^{B_{r},(m)},F_{3}^{B_{r},(m)},...)\overset{\text{law}}{\rightarrow }%
FI^{B_{r}}\text{ as }m\rightarrow \infty .
\end{equation*}%
Motivated by this example, our goal is to characterized in terms of
fragmentation with immigration processes the limiting behavior of%
\begin{equation*}
(m-F_{1}^{(m)},F_{2}^{(m)},F_{3}^{(m)},...)\text{ as }m\rightarrow \infty
\end{equation*}%
for some general fragmentations $F^{(m)}$ where the rates at which particles
split are proportional to a function $\tau $ of their mass. In cases where $%
\tau $ is a power function, this will give the asymptotic behavior of $%
(1-F^{(1)}(\varepsilon \cdot ),F_{2}^{(1)}(\varepsilon \cdot ),...)$ as $%
\varepsilon \rightarrow 0.$

\bigskip

This paper is organized as follows. In the remainder of this section, we
first introduce the fragmentation and fragmentation with immigration
processes we will work with (Subsection 1.1) and then state the main results
on the limiting behavior of $F^{(m)}$ (Subsection 1.2). These results are
proved in Section 2. Sections 3, 4 and 5 are devoted to fragmentations with
a power function $\tau $. Section 3 concerns the behavior near $0$ of such
fragmentations starting from $(1,0,...)$. Section 4 deals with the
asymptotic behavior as $m\rightarrow \infty $ of some CRT representations of
the fragmentations $F^{(m)}$. Section 5 is an application of these results
to a family of fragmentations, namely the ``stable fragmentations'',
introduced by Miermont \cite{mierfmoins},\cite{mierfplus}. Last, Section 6
is an Appendix containing some technical proof and some generalization of
our results to fragmentations with erosion.

\subsection{Fragmentation and fragmentation with immigration processes}

\subsubsection{($\boldsymbol{\tau}$,$\boldsymbol{\nu }$)-fragmentations}

For us, the only distinguishing feature of a particle is its mass, so that
the fragmentation system is characterized at a given time by the decreasing
sequence $s_{1}\geq s_{2}\geq ...\geq 0$ of masses of particles present at
that time. We shall then work in the state space
\begin{equation*}
l_{1}^{\downarrow }:=\left\{ \mathbf{s=(}s_{i})_{i\geq 1}:s_{1}\geq
s_{2}\geq ...\geq 0:\sum\nolimits_{i\geq 1}s_{i}<\infty \right\}
\end{equation*}%
which is equipped with the distance%
\begin{equation*}
d(\mathbf{s,s}^{\prime }):=\sum\nolimits_{i\geq 1}\left| s_{i}-s_{i}^{\prime
}\right| .
\end{equation*}
The \textit{dust }state $(0,0,...)$ is rather denoted by $\mathbf{0}.$
Consider then $(F(t),t\geq 0)$, a c\`{a}dl\`{a}g $l_{1}^{\downarrow }$%
-valued Markov process, and denote by $F^{(m)}$ a version of $F$ starting
from $(m,0,...)$.

\begin{definition}
The process $F$ is called a \textit{fragmentation process} if

$\bullet $ for all $m,t\geq 0$, $\sum\nolimits_{i\geq 1}F_{i}^{(m)}(t)\leq m$

$\bullet $ for all $t_{0}\geq 0$, conditionally on $%
F(t_{0})=(s_{1},s_{2},...)$, $(F(t_{0}+t),t\geq 0)$ is distributed as the
process of the decreasing rearrangements of $F^{(s_{1})}(t)$, $%
F^{(s_{2})}(t) $, ... where the $F^{(s_{i})}$'s are independent versions of $%
F$ starting respectively from $(s_{i},0,0,...),$ $i\geq 1$.
\end{definition}

When $F^{(m)}\overset{\text{law}}{=}mF^{(1)}$ for all $m$, the fragmentation
is usually called \textit{homogeneous.} Such homogeneous processes have been
studied by Bertoin \cite{berthfrag01} and Berestycki \cite{berest02}. In
particular, one knows that when the process is \textit{pure-jump}, its law
is characterized by a so-called \textit{dislocation measure} $\nu $ on
\begin{equation*}
l_{1,\leq 1}^{\downarrow }:=\{\mathbf{s\in }l_{1}^{\downarrow
}:\sum\nolimits_{i\geq 1}s_{i}\leq 1,s_{1}<1\}
\end{equation*}%
that integrates $(1-s_{1})$ and that describes the jumps of the process.
Informally, each mass $s$ will split into masses $ss_{1},ss_{2},...,$ $%
\sum_{i\geq 1}s_{i}\leq 1,$ at rate $\nu ($\textrm{d}$\mathbf{s}).$ We call
such process a $\nu $-\textit{homogeneous fragmentation}. To be more
precise, the papers \cite{berthfrag01},\cite{berest02} give a construction
of the fragmentation based on a Poisson point process $(t_{i},(\mathbf{s(}%
t_{i}),k(t_{i})))_{i\geq 1}$ on $l_{1,\leq 1}^{\downarrow }\times \mathbb{N}$
with intensity measure $\nu \otimes \#$, where $\#$ denotes the counting
measure on $\mathbb{N}$. The construction is so that, at each time $t_{i}$,
the $k(t_{i})$-th mass $F_{k(t_{i})}^{(m)}(t_{i}-)$ splits in masses $%
s_{1}(t_{i})F_{k(t_{i})}^{(m)}(t_{i}-),s_{2}(t_{i})F_{k(t_{i})}^{(m)}(t_{i}-),
$..., the other masses being unchanged. The sequence $F^{(m)}(t_{i})$ is
then the decreasing rearrangement of these new masses and of the unchanged
masses $F_{k}^{(m)}(t_{i}-)$, $k\neq k(t_{i}).$

\textbf{\newline
General setting. }In this paper, we are more generally interested in
pure-jump fragmentation processes where particles with mass $s$ split at
rate $\tau (s)\nu ($\textrm{d}$\mathbf{s})$, where $\tau $ denotes some
continuous strictly positive function on $(0,\infty ).$ When $\nu $ is
finite, this means that each particle with mass $s$ waits an exponential
time with parameter $\tau (s)\nu (l_{1,\leq 1}^{\downarrow })$ before
splitting, and when it splits, it divides into particles with masses $sS_{1}$%
, $sS_{2}$, ..., where $(S_{1},S_{2},...)$ is independent of the splitting
time and is distributed according to $\nu (\cdot )/\nu (l_{1,\leq
1}^{\downarrow })$. When $\nu $ is infinite, the particles split
immediately. In all cases, these models are constructed from homogeneous
fragmentations using time-changes depending on $\tau .$ This is detailed
below. Let us just add here that in the sequel, we will always focus on such
$(\tau ,\nu )$ fragmentations where
\begin{equation}
\begin{array}{l}
\text{- }\tau \text{ is monotone near }0\vspace{0.2cm} \\
\text{- }\nu (\sum\nolimits_{i\geq 1}s_{i}<1)=0,%
\end{array}
\tag{H}  \label{H}
\end{equation}%
the hypothesis on $\nu $ meaning that the fragments do not lose mass within
sudden dislocations.

\textbf{\newline
Construction. }The distribution of each $(\tau ,\nu )$-fragmentation is
constructed through time-changes of a $\nu $-homogeneous fragmentation
starting from $(1,0,...)$ in the following manner (see \cite{Haas1} for
details): let $F^{(1),\text{hom}}$ be a $\nu $-homogeneous fragmentation
starting from $(1,0,...)$ and consider a family $\left( I^{\text{hom}%
}(t),t\geq 0\right) $ of nested random open sets of $(0,1)$ such that $%
F^{(1),\text{hom}}(t)$ is the decreasing sequence of the lengths of interval
components of $I^{\text{hom}}(t)$, for all $t\geq 0.$ One knows (\cite%
{bertsfrag02},\cite{berest02}) that such \textit{interval representation }of
the fragmentation always exists. For $x\in (0,1),$ $t\geq 0,$ call $I_{x}^{%
\text{hom}}(t)$ the connected component of $I^{\text{hom}}(t)$ that contains
$x,$ with the convention $I_{x}^{\text{hom}}(t):=\emptyset $ if $x\notin I^{%
\text{hom}}(t).$ Introduce then the time-changes%
\begin{equation}
T_{x}^{m}(t):=\inf \left\{ u\geq 0:\int_{0}^{u}\frac{\mathrm{d}r}{\tau
(m\left| I_{x}^{\text{hom}}(r)\right| )}>t\right\} ,  \label{4}
\end{equation}%
where $\left| I_{x}^{\text{hom}}(r)\right| $ denotes the length of the
interval $I_{x}^{\text{hom}}(r)$ and, by convention, $\tau (0):=\infty $ and
$\inf \{\emptyset \}:=\infty .$ Clearly, the open sets of $\left( 0,1\right)
$%
\begin{equation*}
I^{\tau }(t):=\bigcup\nolimits_{x\in (0,1)}I_{x}^{\text{hom}}(T_{x}^{m}(t)),%
\text{ \ }t\geq 0,
\end{equation*}%
are nested and we call $F^{(m)}(t)$ the decreasing rearrangement of $m$
times the lengths of the intervals components of $I^{\tau }(t),$ $t\geq 0.$
The process $F^{(m)}$ is then the required fragmentation process starting
from $(m,0,...)$ with splitting rates $\tau (s)\nu ($\textrm{d}$\mathbf{s})$
(Proposition 1, \cite{Haas1}).

\textbf{\newline
Self-similar fragmentations. }When $\tau (s)=s^{\alpha }$ for some $\alpha
\in \mathbb{R}$, the fragmentation is called \textit{self-similar} with
index $\alpha $, since $F^{(m)}\overset{\text{law}}{=}mF^{(1)}(m^{\alpha
}\cdot )$ for all $m>0.$ These self-similar fragmentations processes have
been extensively studied by Bertoin \cite{berthfrag01},\cite{bertsfrag02},%
\cite{bertafrag02}.

\bigskip

\noindent \textbf{Two classical examples. }The \textit{Brownian fragmentation%
} is a self-similar fragmentation process constructed from a normalized
Brownian excursion $e^{(m)}$ with length $m$ as follows: for each $t,$ $%
F^{B_{r},(m)}(t)$ is the decreasing rearrangement of lengths of connected
components of $\{x\in (0,m):2e^{(m)}(x)>t\}$. Equivalently it can be
constructed from the Brownian continuum random tree of Aldous by removing
vertices under height $t$, as explained in the introduction (precise
definition of continuum random trees are given in Section 4). The index of
self-similarity is then $-1/2$ and Bertoin \cite{bertsfrag02} proves that\
the dislocation measure is given by
\begin{equation}
\nu _{B_{r}}\left( s_{1}\in \mathrm{d}x\right) =(2\pi x^{3}\left( 1-x\right)
^{3})^{-1/2}\mathrm{d}x,\text{ }x\in \left[ 1/2,1\right) \text{, and }\nu
_{B_{r}}\left( s_{1}+s_{2}<1\right) =0,  \label{28}
\end{equation}%
this second property meaning that each fragment splits into two pieces when
dislocating.

On the other hand, by logging the Brownian continuum random tree along its
skeleton, Aldous and Pitman \cite{AldousPitman} have introduced a
self-similar fragmentation $F^{AP}$ with index $1/2$ which is transformed by
an exponential time-reversal into the standard additive coalescent. This
\textit{Aldous-Pitman fragmentation} is in some sense dual to the Brownian
one: its dislocation measure is also $\nu _{B_{r}}$ (see \cite{bertsfrag02}).

\textbf{\newline
Loss of mass.} Consider the total mass $M^{(m)}(t)=\sum\nolimits_{i\geq
1}F_{i}^{(m)}(t)$ of macroscopic particles present at time $t$ in a
fragmentation $F^{(m)}$. When the fragmentation rate of small particles is
sufficiently high, some mass may be lost to dust (i.e. a large quantity of
microscopic - or $0$-mass - particles arises in finite time), so that the
mass $M^{(m)}(t)$ decreases to $0$ as $t\rightarrow \infty $. Such
phenomenon does not depend on the initial mass $m$ and happens, for example,
as soon as $\int_{0^{+}}\mathrm{d}x/(x\tau (x))<\infty .$ We refer to \cite%
{Haas1} for some necessary and sufficient condition. An interesting fact is
that the mass $M^{(m)}$ decreases continuously:

\begin{proposition}
\label{Propmassecontinue}The function $t\mapsto M^{(m)}(t)$ is continuous on
$\left[ 0,\infty \right) .$
\end{proposition}

This will be useful for some forthcoming proofs. A proof is given in the
Appendix.

\subsubsection{($\boldsymbol{\tau}$,$\boldsymbol{\nu}$,$\mathbf{I}$%
)-fragmentations with immigration}

Let $\mathcal{I}$ be the set of measures on $l_{1}^{\downarrow }$ that
integrate $(\sum_{j\geq 1}s_{j})\wedge 1$. Two such measures $I,J$ are
considered to be equivalent if their difference $I-J$ puts mass only on $\{%
\mathbf{0\}}$. Implicitly, we always identify a measure with its equivalence
class. In particular, in the following, we will often do the assumption $%
I(l_{1}^{\downarrow })\neq 0,$ which means that $I$ puts mass on some
non-trivial sequences. Endow then $\mathcal{I}$ with the distance%
\begin{equation}
D(I,J)=\sup_{f\in \mathcal{F}}\left| \int_{l_{1}^{\downarrow }}f(\mathbf{s}%
)(I-J)(\text{\textrm{d}}\mathbf{s})\right|  \label{27}
\end{equation}%
where $\mathcal{F}$ is the set of non-negative continuous functions on $%
l_{1}^{\downarrow }$ such that $f(\mathbf{s})\leq (\sum_{j\geq
1}s_{j})\wedge 1$. The function $\mathbf{s}\mapsto (\sum_{j\geq
1}s_{j})\wedge 1$ belongs to $\mathcal{F}$ and therefore $\mathcal{I}$ is
closed. It is called the set of \textit{immigration measures}.

\begin{definition}
Let $\left( \left( r_{i},\mathbf{u}^{i}\right) ,i\geq 1\right) $ be a
Poisson point process (PPP) with intensity $I\in \mathcal{I}$ and,
conditionally on this PPP, let $F^{(u_{j}^{i})},$ $i,j\geq 1,$ be
independent $\left( \tau ,\nu \right) $ fragmentations starting respectively
from $u_{j}^{i},$ $i,j\geq 1$. Then consider for each $t\geq 0,$ the
decreasing rearrangement%
\begin{equation*}
FI(t):=\left\{ F_{k}^{(u_{j}^{i})}(t-r_{i}),r_{i}\leq t,j,k\geq 1\right\}
^{\downarrow }\in l_{1}^{\downarrow }\text{.}
\end{equation*}%
The process $FI$ is called a fragmentation with immigration process with
parameters $\left( \tau ,\nu ,I\right) $.

When there is no fragmentation ($\nu (l_{1,\leq 1}^{\downarrow })=0$), we
rather call such process a pure immigration process with parameter $I$ and
we denote it by $(I(t),t\geq 0).$
\end{definition}

This means that at time $r_{i}$, particles with masses $u_{1}^{i}$, $%
u_{2}^{i}$, ... immigrate and then start to fragment independently of each
other (conditionally on their masses), according to a $\left( \tau ,\nu
\right) $ fragmentation. The initial state is $\mathbf{0}$. Note that the
total mass of immigrants until time $t$%
\begin{equation}
\sigma _{I}(t):=\sum\nolimits_{r_{i}\leq t,j\geq 1}u_{j}^{i}  \label{5}
\end{equation}%
is a.s. finite and therefore that the decreasing rearrangement $FI(t)$
indeed exists and is in $\in l_{1}^{\downarrow }$. The process $\sigma _{I}$
is a \textit{subordinator,} i.e. an increasing L\'{e}vy process. We refer to
the lecture \cite{Bertoin3}, for backgrounds on subordinators. In
particular, we recall that a subordinator $\sigma $ is characterized by its
Laplace exponent, which is a function $\phi _{\sigma }$ such that $E[\exp
(-q\sigma (t))]=\exp (-t\phi _{\sigma }(q))$, for all $q,t\geq 0$.

Note also that $FI$ is c\`{a}dl\`{a}g, since the $F^{(u_{j}^{i})}$ are c\`{a}%
dl\`{a}g, since dominated convergence applies and since, clearly, the
following result holds.

\begin{lemma}
\label{Lemmadecresase}For all integers $1\leq n\leq \infty $, let $%
x^{n}=(x_{i}^{n},i\geq 1)$ be a sequence of non-negative real numbers such
that $\sum\nolimits_{i\geq 1}x_{i}^{n}<\infty $ and let $x^{n}{}^{\downarrow
}$ denotes its decreasing rearrangement. If $\sum\nolimits_{i\geq 1}\left|
x_{i}^{n}-x_{i}^{\infty }\right| \rightarrow 0$, then $\sum\nolimits_{i\geq
1}|x_{i}^{n\downarrow }-x_{i}^{\infty \downarrow }|\rightarrow 0$, i.e. $%
x^{n}{}^{\downarrow }\rightarrow x^{\infty \downarrow }$ in $%
l_{1}^{\downarrow }$.
\end{lemma}

Equilibrium for such fragmentation with immigration processes has been
studied in \cite{HaasImmig04} in a slightly less general context.

\subsection{Main results: asymptotics of $F^{(m)}$}

From now on, we suppose that $\nu (l_{1,\leq 1}^{\downarrow })\neq 0$.
Introduce then for all $m\geq 0$, the measure $\nu _{m}\in \mathcal{I}$
defined for all non-negative measurable functions $f$ on $l_{1}^{\downarrow
} $ by
\begin{equation*}
\int_{l_{1}^{\downarrow }}f(\mathbf{s})\nu _{m}(\text{\textrm{d}}\mathbf{s}%
):=\int_{l_{1,\leq 1}^{\downarrow }}f(s_{2}m,s_{3}m,...)\nu (\text{\textrm{d}%
}\mathbf{s}).
\end{equation*}%
Set also%
\begin{equation*}
\varphi _{\nu }(m):=\left( \nu \left( s_{1}<1-m^{-1}\right) \right)
^{-1}=\langle \nu _{m},\mathbf{1}_{\{\sum_{i\geq 1}s_{i}>1\}}\rangle ^{-1}
\end{equation*}%
which is finite for $m$ large enough and converges to $0$ as $m\rightarrow
\infty $ when $\nu (l_{1,\leq 1}^{\downarrow })=\infty .$

We are now ready to state our main result. We remind that the distance on $%
\mathcal{I}$ is defined by $\left( \ref{27}\right) $. Also, the set of c\`{a}%
dl\`{a}g paths in $\mathbb{R}^{+}\mathbb{\times }l_{1}^{\downarrow }$ is
endowed with the \textit{Skorohod topology.}

\begin{theorem}
\label{Thprincipal}Let $F$ be a $(\tau ,\nu )$ fragmentation and suppose
that $\tau (m)\nu _{m}\rightarrow I$, $I(l_{1}^{\downarrow })\neq 0$, as $%
m\rightarrow \infty $. Then,%
\begin{equation*}
\left( m-F_{1}^{(m)},(F_{2}^{(m)},F_{3}^{(m)},...)\right) \overset{\mathrm{%
law}}{\rightarrow }(\sigma _{I},FI)\text{ as }m\rightarrow \infty
\end{equation*}%
where $FI$ is a fragmentation with immigration with parameters $\left( \tau
,\nu ,I\right) $, starting from $\mathbf{0}$ and $\sigma _{I}$ is the
process $\left( \ref{5}\right) $ corresponding to the total mass of
particles that have immigrated until time $t$, $t\geq 0$.
\end{theorem}

In some sense, letting $m\rightarrow \infty $ in $F^{(m)}$ creates an
infinite amount of mass that regularly injects into the system some groups
of finite masses which then undergo fragmentation. A similar phenomenon has
been observed in the study of some different processes conditioned on
survival (see e.g. \cite{DuquesneImmigration},\cite{EthridgeWilliams},\cite%
{Evans},\cite{Lyons Pemantle Peres}).\bigskip

\noindent \textbf{Example.} Recall the characterization (\ref{28}) of the
Brownian dislocation measure $\nu _{B_{r}}$. Clearly, $m^{-1/2}\nu
_{B_{r},m}\rightarrow I_{B_{r}}$ where the measure $I_{B_{r}}$ is defined by
\begin{equation}
I_{B_{r}}(s_{1}\in \mathrm{d}x)=(2\pi x^{3})^{-1/2}\mathrm{d}x,\text{ }x>0,%
\text{ and }I_{B_{r}}(s_{2}>0)=0.  \label{29}
\end{equation}%
So the previous theorem applies to the Brownian fragmentation and the
fragmentation with immigration appearing in the limit has parameters $\left(
\tau :x\mapsto x^{-1/2},\nu _{B_{r}},I_{B_{r}}\right) $. The L\'{e}vy
measure of the subordinator $\sigma _{I_{B_{r}}}$ is simply $%
I_{B_{r}}(s_{1}\in \mathrm{d}x).$ Informally, this corresponds to the
convergence, mentioned in the introduction, of the Brownian CRT to a tree
with a spine on which are branched rescaled Brownian CRTs. This tree with a
spine codes (see Section 4 for precise statements) the above $\left( \tau
:x\mapsto x^{-1/2},\nu _{B_{r}},I_{B_{r}}\right) $ fragmentation with
immigration.

\vspace{0.2cm}Other examples are given in Section 5.1.\bigskip

The assumption on the convergence of $\tau (m)\nu _{m}$ may seem demanding
and, clearly, is not always satisfied. A moment of thought, using
test-functions of type $f_{a}(\mathbf{s})=\mathbf{1}_{\{\sum_{i\geq
1}s_{i}>a\}}$, $a>0$, leads to the following result.

\begin{lemma}
\label{Lemma 6}Suppose that $\tau (m)\nu _{m}$ converges to some measure $I$%
, $I(l_{1}^{\downarrow })\neq 0$, as $m\rightarrow \infty $. Then both $\tau
$ and $\varphi _{\nu }$ vary regularly at $\infty $ with some index $-\gamma
_{\nu },$ $\gamma _{\nu }\in (0,1)$ and $\tau (m)\sim C\varphi _{\nu }(m)$
as $m\rightarrow \infty $, for some constant $C>0$. As a consequence, the
limit $I$ is $\gamma _{\nu }$\textit{-self-similar}, that is
\begin{equation*}
\int_{l_{1}^{\downarrow }}f(as_{1},as_{2},...)I\left( \mathrm{d}\mathbf{s}%
\right) =a^{\gamma _{\nu }}\int_{l_{1}^{\downarrow }}f(\mathbf{s})I\left(
\mathrm{d}\mathbf{s}\right) \text{ for all }a>0\text{, }f\in \mathcal{F}%
\text{,}
\end{equation*}%
which in turn implies that $\sigma _{I}$ is a stable subordinator with index
$\gamma _{\nu }$ and Laplace exponent $C\Gamma (1-\gamma _{\nu })q^{\gamma
_{\nu }}$, $q\geq 0$.
\end{lemma}

Hence Theorem \ref{Thprincipal} applies to measures $\nu $ such that $%
\varphi _{\nu }(m)\nu _{m}$ converges, coupled together with functions $\tau
$ whose behavior at $\infty $ is proportional to that of $\varphi _{\nu }$.
Note in particular that the speed of fragmentation of small particles plays
no role in the existence of a limit.

\bigskip

Remark then that it is possible to construct from any $\gamma $-self-similar
immigration measure $I,$ $I(l_{1}^{\downarrow })\neq 0,$ $\gamma \in (0,1),$
some dislocation measures $\nu $ such that $\varphi _{\nu }(m)\nu _{m}$
converge\footnote{%
For example, define $\nu $ by $\int_{\mathcal{S}^{\downarrow }}f(\mathbf{s}%
)\nu (\mathrm{d}\mathbf{s}):=\int_{l_{1}^{\downarrow
}}f(1-\sum\nolimits_{j\geq 1}s_{j},s_{1},s_{2},...)\mathbf{1}_{\{s_{1}\leq
1-\sum_{j\geq 1}s_{j}\}}I(\mathrm{d}\mathbf{s}).$ Clearly, $\nu (\sum_{j\geq
1}s_{j}\neq 1)=0$, $\nu $ integrates $(1-s_{1})$ and $m^{-\gamma }\nu
_{m}\rightarrow I$.
\par
{}} to $I$, which gives a large class of measures $\nu $ to which Theorem %
\ref{Thprincipal} applies. Also, note that when the fragmentation is binary
(i.e. when $\nu (s_{1}+s_{2}<1)=0$), the convergence of $\varphi _{\nu
}(m)\nu _{m}$ holds as soon as $\varphi _{\nu }$ varies regularly at $\infty
$ with some index in $(-1,0).$

\bigskip

For functions $\tau $ such that $(\varphi _{\nu }/\tau )(m)$ converges to $0$
or $\infty $, a first computation shows that, provided $\varphi _{\nu
}(m)\nu _{m}$ converges and $\tau $ varies regularly at $\infty $:

- $(\varphi _{\nu }/\tau )(m)\rightarrow 0\Rightarrow
(m-F_{1}^{(m)}(1),F_{2}^{(m)}(1))\overset{\mathrm{law}}{\rightarrow }(\infty
,\infty )$

- $(\varphi _{\nu }/\tau )(m)\rightarrow \infty \Rightarrow
(m-F_{1}^{(m)}(1),(F_{2}^{(m)}(1),F_{3}^{(m)}(1),...))\overset{\mathrm{law}}{%
\rightarrow }(0,\mathbf{0})$.

\bigskip

One way to avoid these trivial limits is to consider the process $F^{(m)}$
up to a time change:

\begin{theorem}
\label{Proposition_infinie} Suppose that $\tau $ varies regularly at $\infty
$, and that $\varphi _{\nu }(m)\nu _{m}\rightarrow I$ as $m\rightarrow
\infty $.

\noindent $\mathrm{(i)}$ If $(\varphi _{\nu }/\tau )(m)\rightarrow 0$, then,
as $m\rightarrow \infty $,\vspace{-0.15cm}
\begin{equation*}
((m-F_{1}^{(m)}\left( (\varphi _{\nu }/\tau )(m)\cdot \right)
),F_{2}^{(m)}\left( (\varphi _{\nu }/\tau )(m)\cdot \right)
,F_{3}^{(m)}\left( (\varphi _{\nu }/\tau )(m)\cdot \right) ,...)\overset{%
\mathrm{law}}{\rightarrow }(\sigma _{I},(I(t),t\geq 0)),\vspace{-0.15cm}
\end{equation*}%
where $(I(t),t\geq 0)$ is a pure immigration process with parameter $I$.

\vspace{0.2cm}\noindent $\mathrm{(ii)}$ If $(\varphi _{\nu }/\tau
)(m)\rightarrow \infty $ and the fragmentation loses mass to dust, then the
following finite-dimensional convergence holds as $m\rightarrow \infty ,$%
\vspace{-0.15cm}
\begin{equation*}
((m-F_{1}^{(m)}\left( (\varphi _{\nu }/\tau )(m)\cdot \right)
),F_{2}^{(m)}\left( (\varphi _{\nu }/\tau )(m)\cdot \right)
,F_{3}^{(m)}\left( (\varphi _{\nu }/\tau )(m)\cdot \right) ,...)\overset{%
\mathrm{law}}{\underset{\mathrm{f.d.}}{\rightarrow }}(\sigma _{I},\mathbf{0}%
).
\end{equation*}
\end{theorem}

The assertion (ii) is not valid when the fragmentation does not lose mass,
since the quantity $m-\sum\nolimits_{i\geq 1}F_{i}^{(m)}\left( (\varphi
_{\nu }/\tau )(m)\right) $ is then equal to $0$ and so cannot converge to $%
\sigma _{I}(1).$ However, a result similar to that stated in (ii) holds for
fragmentations that do not lose mass, provided that the distance $d$ is
replaced by the distance of uniform convergence on $l_{1}^{\downarrow }.$
Also, the reason why the limit in this statement (ii) holds only in the
finite dimensional sense and not with respect to the Skorohod topology, is,
informally, that the functional limit of $(F_{2}^{(m)}\left( (\varphi _{\nu
}/\tau )(m)\cdot \right) ,F_{3}^{(m)}\left( (\varphi _{\nu }/\tau )(m)\cdot
\right) ,...)$ cannot be c\`{a}dl\`{a}g.

Another remark is that under the assumptions of the first statement, the
processes $F_{i}^{(m)}((\varphi _{\nu }/\tau )(m)\cdot ),$ $i\geq 2,$
although not increasing, converge as $m\rightarrow \infty $ to some
increasing processes. In particular, $F_{2}^{(m)}((\varphi _{\nu }/\tau
)(m)\cdot )$ converges to $(\Delta _{1}^{a}(t),t\geq 0)$ where $\Delta
_{1}^{a}(t)$ is the largest jump before time $t$ of some stable subordinator
with Laplace exponent $a\Gamma (1-\gamma _{\nu })q^{\gamma _{\nu }}$, $q\geq
0$, and $a=\lim_{m\rightarrow \infty }\varphi _{\nu }(m)\nu (s_{2}>m^{-1})$
(this limit exists, although $\mathbf{s}\mapsto \mathbf{1}%
_{\{s_{1}>1\}}\notin \mathcal{F},$ because $I(s_{1}\in \mathrm{d}x)$ is
absolutely continuous, as a consequence of the self-similarity). In case $%
\nu $ is binary, one more precisely has:

\begin{corollary}
\label{Corosub}Suppose that $\nu $ is binary and suppose that $\varphi _{\nu
}$ varies regularly at $\infty $ with some index $-\gamma _{\nu },\gamma
_{\nu }\in (0,1)$. Then, if $\alpha >-\gamma _{\nu }$,
\begin{equation*}
((m-F_{1}^{(m)}\left( (\varphi _{\nu }/\tau )(m)\cdot \right)
),F_{2}^{(m)}\left( (\varphi _{\nu }/\tau )(m)\cdot \right)
,F_{3}^{(m)}\left( (\varphi _{\nu }/\tau )(m)\cdot \right) ,...)\overset{%
\mathrm{law}}{\rightarrow }(\sigma ,\Delta _{1},\Delta _{2},...)
\end{equation*}%
where $\sigma $ is a stable subordinator with Laplace exponent $\Gamma
(1-\gamma _{\nu })q^{\gamma _{\nu }}$ and $(\Delta _{1}(t),\Delta
_{2}(t),...)$ the decreasing sequence of its jumps before time $t,$ $t\geq 0$%
.
\end{corollary}

\noindent \textbf{Example. }This can be applied to the Aldous-Pitman
fragmentation, since $\varphi _{\nu _{B_{r}}}(m)\sim \pi ^{1/2}(2m)^{-1/2}$.
We get that
\begin{equation*}
((m-F_{1}^{AP,(m)}(m^{-1}\cdot )),F_{2}^{AP,(m)}(m^{-1}\cdot
),F_{3}^{AP,(m)}(m^{-1}\cdot ),...)\overset{\mathrm{law}}{\rightarrow }%
(\sigma _{AP},\Delta _{1}^{AP},\Delta _{2}^{AP},...)
\end{equation*}%
where $\sigma _{AP}$ is a stable subordinator with Laplace exponent $\sqrt{2}%
q^{1/2}$ and $(\Delta _{1}^{AP}(t),\Delta _{2}^{AP}(t),...)$ the decreasing
sequence of its jumps before time $t$, $t\geq 0$. Aldous and Pitman \cite%
{AldousPitman}, Corollary 13, obtained this result by studying size-biased
permutations of their fragmentation.

\vspace{0.2cm}Other explicit (and non-binary) examples are studied in
Section 5.2.

\section{Proofs}

The main lines of the proofs of Theorems \ref{Thprincipal} and \ref%
{Proposition_infinie} are quite similar. We first give a detailed proof of
Theorem \ref{Thprincipal} and then explain how it can be adapted to prove
Theorem \ref{Proposition_infinie}. We will need the following classical
result on Skorohod convergence (see Proposition 3.6.5, \cite{Ethier Kurtz}).

\begin{lemma}
\label{LemmaEthr}Consider a metric space $(E,d_{E})$ and let $f_{n},f$ be c%
\`{a}dl\`{a}g paths with values in $E.$ Then $f_{n}\rightarrow f$ with
respect to the Skorohod topology if and only if the three following
assertions are satisfied for all sequences $t_{n}\rightarrow t$, $%
t_{n},t\geq 0:$

$\mathrm{(a)}$ $\min
(d_{E}(f_{n}(t_{n}),f(t)),d_{E}(f_{n}(t_{n}),f(t-)))\rightarrow 0$

$\mathrm{(b)}$ $d_{E}(f_{n}(t_{n}),f(t))\rightarrow 0$ $\Rightarrow $ $%
d_{E}(f_{n}(s_{n}),f(t))\rightarrow 0$ for all sequences $s_{n}\rightarrow
t, $ $s_{n}\geq t_{n}$

$\mathrm{(c)}$ $d_{E}(f_{n}(t_{n}),f(t-))\rightarrow 0$ $\Rightarrow $ $%
d_{E}(f_{n}(s_{n}),f(t-))\rightarrow 0$ for all sequences $s_{n}\rightarrow
t,$ $s_{n}\leq t_{n}$
\end{lemma}

\subsection{Proof of Theorem \ref{Thprincipal}}

In this section it is supposed that $\tau (m)\nu _{m}\rightarrow I$, $%
I(l_{1}^{\downarrow })\neq 0$, as $m\rightarrow \infty .$ Our goal is then
to prove Theorem \ref{Thprincipal}, which is a corollary of the forthcoming
Lemma \ref{LemmaJoint}. In order to state and prove this lemma, we first
introduce some notations and give some heuristic geometric description of
what is happening. There is no loss of generality in supposing that the $%
(\tau ,\nu )$ fragmentations $F^{(m)}$, $m\geq 0$, are constructed from the
same $\nu $-homogeneous one, which is done in the following.

\subsubsection{Heuristic description}

We first give a geometric description of the fragmentation $F^{(m)}$, which
may be viewed as a baseline $\mathcal{B}=\left[ 0,\infty \right) $ on which
fragmentation processes are attached.

Let $\Lambda ^{(m)}$ be the process obtained by following at each
dislocation the largest sub-fragment. According to the Poissonian
construction of homogeneous fragmentation processes and the time-change
between $\nu $-homogeneous and $(\tau ,\nu )$-fragmentations (see Section
1.1.1), the process $\Lambda ^{(m)}$ is constructed from some Poisson point
process $(\left( t_{i},\mathbf{s}^{i}\right) ,i\geq 1)$ (independent of $m$)
with intensity measure $\nu $ as follows: if $\xi $ denotes the subordinator
defined by
\begin{equation}
\xi (t):=\sum\nolimits_{t_{i}\leq t}(-\log (s_{1}^{i})),\text{ }t\geq 0,
\label{31}
\end{equation}%
and $\rho ^{(m)}$ the time change
\begin{equation}
\rho ^{(m)}(t):=\inf \left\{ u:\int_{0}^{u}\mathrm{d}r/\tau (m\exp (-\xi
(r)))>t\right\} ,  \label{16}
\end{equation}%
then%
\begin{equation}
\Lambda ^{(m)}(t)=m\exp (-\xi (\rho ^{(m)}(t)))\text{, }t\geq 0\text{.}
\label{2}
\end{equation}%
The set of jump times of $\Lambda ^{(m)}$ is $\{t_{i}^{m}:=\rho
^{-(m)}(t_{i}),i\geq 1\}$.

The evolution of the fragmentation $F^{(m)}$ then relies on the point
process $((t_{i}^{m},\mathbf{s}^{i}),i\geq \nolinebreak 1)$: at time $%
t_{i}^{m}$, the fragment with mass $\Lambda ^{(m)}(t_{i}^{m}-)$ splits to
give a fragment with mass $\Lambda ^{(m)}(t_{i}^{m})=\Lambda
^{(m)}(t_{i}^{m}-)s_{1}^{i}$ and smaller fragments with masses $\Lambda
^{(m)}(t_{i}^{m}-)s_{j}^{i}$, $j\geq 2.$ For $j\geq 2$, call $F^{(\Lambda
^{(m)}(t_{i}^{m}-)s_{j}^{i})}$ the fragmentation describing the evolution of
the mass $\Lambda ^{(m)}(t_{i}^{m}-)s_{j}^{i}$ and consider that it is
branched at height $t_{i}^{m}$ on the baseline $\mathcal{B}$. Then the
process $F^{(m)}$ is obtained by considering for each $t\geq 0$ all
fragmentations branched at height $t_{i}^{m}\leq t$ and by ordering in the
decreasing order the terms of sequences $F^{(\Lambda
^{(m)}(t_{i}^{m}-)s_{j}^{i})}(t-t_{i}^{m}),$ $t_{i}^{m}\leq t,$ $j\geq 2,$
and $\Lambda ^{(m)}(t)$. In some sense, there is then a tree structure under
this baseline with ``fragmentation'' leaves. This will be discussed in
Section 4.

Similarly, a $(\tau ,\nu ,I)$ fragmentation with immigration $FI$ can be
viewed as the baseline $\mathcal{B}$ with fragmentations leaves $%
F^{(u_{j}^{i})}$, $j\geq 1,$ attached at time $r_{i}$, where $((r_{i},%
\mathbf{u}^{i}),i\geq 1)$ is a Poisson point process with intensity $I$ and $%
F^{(u_{j}^{i})},$ $i,j\geq 1,$ some $(\tau ,\nu )$ fragmentations starting
respectively from $u_{j}^{i},$ $i,j\geq 1,$ that are independent
conditionally on $((r_{i},\mathbf{u}^{i}),i\geq 1).$

Now, to see the connection between these descriptions and the result we want
to prove on the convergence of $(F_{2}^{(m)},F_{3}^{(m)},...)$ to $FI$, note
that the processes $\Lambda ^{(m)}$ and $F_{1}^{(m)},$ although different,
coincide at least when $\Lambda ^{(m)}(t)\geq m/2,$ since $\Lambda ^{(m)}(t)$
is then the largest fragment of $F^{(m)}(t)$. Fix $t_{0}<\infty $.\ It is
easily seen that under the assumption $\tau (m)\nu _{m}\rightarrow I$ (which
in particular implies that $\tau (m)\rightarrow 0$ as $m\rightarrow \infty $%
), a.s. $\rho ^{(m)}(t_{0})\rightarrow 0$ as $m\rightarrow \infty $, which
in turn implies that for large $m$'s and all $t\leq t_{0}$, $\Lambda
^{(m)}(t)\geq m/2$, and therefore $\Lambda ^{(m)}(t)=F_{1}^{(m)}(t)$. In
particular $(F_{2}^{(m)}(t),F_{3}^{(m)}(t),...)$ is then the decreasing
rearrangement of the terms of sequences $F^{(\Lambda
^{(m)}(t_{i}^{m}-)s_{j}^{i})}(t-t_{i}^{m}),$ $t_{i}^{m}\leq t,$ $j\geq 2.$

Hence, informally, one may expect that the process $%
(F_{2}^{(m)},F_{3}^{(m)},..)$ converges in law to $FI$ as soon as $(\Lambda
^{(m)}(t_{i}^{m}-)s_{j}^{i},i\geq 1,j\geq 2)$ converges to $%
(u_{j}^{i},i,j\geq 1),$ and $(t_{i}^{m},i\geq 1)$ to $(r_{i},i\geq 1)$. The
statement of these convergences is made rigorous in the forthcoming Lemma %
\ref{Lemmamesures}, which is then used to prove the required Lemma \ref%
{LemmaJoint}.

\subsubsection{Convergence of the point processes}

Consider the set $\left[ 0,\infty \right) \times l_{1}^{\downarrow }$
endowed with the product topology (which makes it Polish) and introduce the
set $\mathcal{R}_{\left[ 0,\infty \right) \times l_{1}^{\downarrow }}$ of
Radon point measures on $\left[ 0,\infty \right) \times l_{1}^{\downarrow }$
that integrate $\mathbf{1}_{\left\{ t\leq t_{0}\right\} }\times
\sum\nolimits_{j\geq 1}s_{j},$ for all $t_{0}\geq 0.$ Two such measures are
considered to be equivalent if their difference puts mass only on $\left[
0,\infty \right) \times \{\mathbf{0}\}$. Again, we shall implicitly identify
a measure with its equivalence class. Introduce then $\mathcal{F}_{\left[
0,\infty \right) \times l_{1}^{\downarrow }}$, the set of $\mathbb{R}^{+}$%
-valued continuous functions $f$ on $\left[ 0,\infty \right) \times
l_{1}^{\downarrow }$ such that $f(t,\mathbf{s})\leq \mathbf{1}_{\left\{
t\leq t_{0}\right\} }\sum\nolimits_{j\geq 1}s_{j}$ for some $t_{0}\geq 0$
(we shall denote by $t_{0}^{f}$ such $t_{0}$'s) and equip $\mathcal{R}_{%
\left[ 0,\infty \right) \times l_{1}^{\downarrow }}$ with the topology
induced by the convergence $\mu _{n}\rightarrow \mu \Leftrightarrow \langle
\mu _{n},f\rangle \rightarrow \langle \mu ,f\rangle $ for all $f\in \mathcal{%
F}_{\left[ 0,\infty \right) \times l_{1}^{\downarrow }}$. With respect to
this topology, one has

\begin{lemma}
\label{Lemmamesures}$\sum\nolimits_{i\geq 1}\delta _{(t_{i}^{m},(\Lambda
^{(m)}(t_{i}^{m}-)s_{j}^{i})_{j\geq 2})}\overset{\mathrm{law}}{\rightarrow }%
\sum\nolimits_{i\geq 1}\delta _{(r_{i},\mathbf{u}^{i})}$ as $m\rightarrow
\infty .$
\end{lemma}

\noindent \textbf{Proof. }We first point out that $\mu
_{m}:=\sum\nolimits_{i\geq 1}\delta _{(t_{i}/\tau (m),(ms_{j}^{i})_{j\geq
2})}$ converges in distribution to $\mu :=\sum\nolimits_{i\geq 1}\delta
_{(r_{i},\mathbf{u}^{i})}$. Indeed, both measures belong to $\mathcal{R}_{%
\left[ 0,\infty \right) \times l_{1}^{\downarrow }}$ and according to
Theorems 4.2 and 4.9 of Kallenberg \cite{Kallenberg}, the convergence in
distribution of $\mu _{m}$ to $\mu $ is equivalent to the convergence of all
Laplace transforms $E\left[ \exp (-\langle \mu _{m},f\rangle )\right] $ to $E%
\left[ \exp (-\langle \mu ,f\rangle )\right] ,$ $f\in \mathcal{F}_{\left[
0,\infty \right) \times l_{1}^{\downarrow }}$, which is easily checked: fix
such function $f$ and apply Campbell formula (see e.g. \cite{kingman93}) to
the Poisson point processes $((t_{i},\mathbf{s}^{i}),i\geq 1)$ to obtain%
\begin{equation*}
E\left[ \exp (-\langle \mu _{m},f\rangle )\right] =\exp \left( -\tau
(m)\int_{\left[ 0,t_{0}^{f}\right] \times l_{1}^{\downarrow }}(1-\exp (-f(u,%
\mathbf{s})))(\text{\textrm{d}}u\otimes \nu _{m}(\text{\textrm{d}}\mathbf{s}%
))\right) .
\end{equation*}%
Clearly, the function $F_{f}:\mathbf{s\mapsto }\int_{0}^{t_{0}^{f}}(1-\exp
(-f(u,\mathbf{s})))$\textrm{d}$u$ is continuous and bounded by $%
t_{0}^{f}((\sum\nolimits_{j\geq 1}s_{j})\wedge 1).$ Therefore $\langle \tau
(m)\nu _{m},F_{f}\rangle \rightarrow \langle I,F_{f}\rangle $, which in turn
implies that \linebreak $E\left[ \exp (-\langle \mu _{m},f\rangle )\right]
\rightarrow E\left[ \exp (-\langle \mu ,f\rangle )\right] $. Hence $\mu _{m}%
\overset{\mathrm{law}}{\rightarrow }\mu $.

Then, using Skorohod's representation theorem (our set of point measures is
Polish, see e.g. Appendix A7 of Kallenberg \cite{Kallenberg}), one may
suppose that $\mu _{m}\rightarrow \mu $ a.s. To simplify, we work in the
rest of the proof with the representation of the measure $\mu $ (resp. $\mu
_{m}$, $m\geq 0$) that does not put mass on $\left[ 0,\infty \right) \times
\{\mathbf{0}\}$.\ We then call $\sigma ^{m}$ a (random) permutation such
that $t_{\sigma ^{m}(i)}/\tau (m)\rightarrow r_{i}$ and $(ms_{j}^{\sigma
^{m}(i)})_{j\geq 2})\rightarrow \mathbf{u}^{i}$, $\forall i\geq 1$, a.s.
This leads us to the a.s. pointwise convergence

\begin{eqnarray}
r_{i}^{m} &:&=t_{\sigma ^{m}(i)}^{m}=\rho ^{-(m)}(t_{\sigma
^{m}(i)})\rightarrow r_{i}  \label{3} \\
\mathbf{z}^{i,m} &:&=(ms_{j}^{\sigma ^{m}(i)})_{j\geq 2}\rightarrow \mathbf{u%
}^{i}  \notag \\
\mathbf{u}^{i,m} &:&=(\Lambda ^{(m)}(t_{\sigma ^{m}(i)}^{m}-)s_{j}^{\sigma
^{m}(i)})_{j\geq 2}\rightarrow \mathbf{u}^{i}.  \notag
\end{eqnarray}%
Indeed, as noticed in Lemma \ref{Lemma 6}, the assumption $\tau (m)\nu
_{m}\rightarrow I,$ $I(l_{1}^{\downarrow })\neq 0$, implies that $\tau $
varies regularly at $\infty $ with index $-\gamma _{\nu }\in (0,1)$. In
particular $\tau (m)\rightarrow 0$ and then $t_{\sigma ^{m}(i)}\rightarrow 0$%
. This implies that $\Lambda ^{(m)}(t_{\sigma ^{m}(i)}^{m}-)=m\exp (-\xi
(t_{\sigma ^{m}(i)}-))\sim m$ and then that $\mathbf{u}^{i,m}\rightarrow
\mathbf{u}^{i}.$ Next, because of the regular variation of $\tau $, one
knows (see Potter's Theorem, Th.1.5.6 \cite{bgt}) that there exists for all $%
A>1,\varepsilon >0,$ some constant $M(A,\varepsilon )\geq 0$ such that%
\begin{equation}
A^{-1}\exp ((-\gamma _{\nu }-\varepsilon )\xi (r))\leq \frac{\tau (m)}{\tau
(m\exp (-\xi (r)))}\leq A\exp ((-\gamma _{\nu }+\varepsilon )\xi (r))\text{,}
\label{10}
\end{equation}%
for all $m,r$ such that $m\exp (-\xi (r))\geq M(A,\varepsilon ).$ This
implies that
\begin{equation*}
\tau (m)\int_{0}^{t_{\sigma ^{m}(i)}}\mathrm{d}r/\tau (m\exp (-\xi (r)))%
\underset{\infty }{\sim }t_{\sigma ^{m}(i)}
\end{equation*}%
and therefore that $r_{i}^{m}=\rho ^{-(m)}(t_{\sigma ^{m}(i)})\rightarrow
r_{i}$ as $m\rightarrow \infty $.

It remains to prove that (a.s.) $\sum\nolimits_{i\geq 1}\delta _{(r_{i}^{m},%
\mathbf{u}^{i,m})}\rightarrow \sum\nolimits_{i\geq 1}\delta _{(r_{i},\mathbf{%
u}^{i})}$. So fix $f\in \mathcal{F}_{\left[ 0,\infty \right) \times
l_{1}^{\downarrow }}$ and choose some $t_{0}^{f}$ and $C>1$ such that $%
t_{0}^{f}C\notin \{r_{i},i\geq 1\}$. Then fix $\eta >0$ and take $i_{0}$
such that $\sum\nolimits_{i>i_{0},r_{i}\leq t_{0}^{f}C}\sum\nolimits_{j\geq
1}u_{j}^{i}<\eta .$\ Since $\mu _{m}\rightarrow \mu $ and $t_{0}^{f}C\neq
r_{i},$ $i\geq 1$, one has%
\begin{equation*}
\sum\nolimits_{i>i_{0},(t_{\sigma ^{m}(i)}/\tau (m))\leq
t_{0}^{f}C}\sum\nolimits_{j\geq 2}(ms_{j}^{\sigma ^{m}(i)})<\eta
\end{equation*}%
for $m$ large enough. Next, we claim that there exists some $m_{0}$ such
that for all $m\geq m_{0}$ and $t\geq 0$, $\rho ^{-(m)}(t)\leq t_{0}^{f}$
leads to $(t/\tau (m))\leq t_{0}^{f}C$, at least if $C$ has been chosen
large enough. Indeed, for $m$ large enough, the left hand side of $\left( %
\ref{10}\right) $ is valid uniformly in $t$, $\forall t\leq 1$. Taking $C$
larger if necessary, we get that $(t/\tau (m))C^{-1}\leq \rho ^{-(m)}(t)$, $%
\forall t\leq 1$, hence that $\rho ^{-(m)}(t)\leq t_{0}^{f}$ implies $%
(t/\tau (m))\leq Ct_{0}^{f}$, $\forall t\leq 1$. On the other hand, still
for $m$ large enough, $\rho ^{-(m)}(1)>t_{0}^{f}$ (since $\tau
(m)\rightarrow 0$), and a fortiori $\rho ^{-(m)}(t)>t_{0}^{f}$ for all $%
t\geq 1$. This leads to the claim. So, for $m$ large enough, $r_{i}^{m}\leq
t_{0}^{f}\ $implies $(t_{\sigma ^{m}(i)}/\tau (m))\leq t_{0}^{f}C$, and
therefore%
\begin{equation*}
\sum\nolimits_{i>i_{0},r_{i}^{m}\leq t_{0}^{f}}\sum\nolimits_{j\geq
2}(ms_{j}^{\sigma ^{m}(i)})<\eta .
\end{equation*}%
Consequently (using that $u_{j}^{i,m}\leq ms_{j+1}^{\sigma ^{m}(i)}$, $j\geq
1$),%
\begin{equation*}
\left| \sum\nolimits_{i\geq 1}(f(r_{i}^{m},\mathbf{u}^{i,m})-f(r_{i},\mathbf{%
u}^{i})\right| \leq \sum\nolimits_{i\leq i_{0}}\left| f(r_{i}^{m},\mathbf{u}%
^{i,m})-f(r_{i},\mathbf{u}^{i})\right| +2\eta .
\end{equation*}%
At last, using first the pointwise convergence (\ref{3}) for the finite sum
on $i\leq i_{0}$ and then letting $\eta \rightarrow 0$, one obtains the
required $\sum\nolimits_{i\geq 1}f(r_{i}^{m},\mathbf{u}^{i,m})\rightarrow
\sum\nolimits_{i\geq 1}f(r_{i},\mathbf{u}^{i}).$ Let us also point out that,
exactly in the same way, one obtains $\sum\nolimits_{i\geq 1}\delta
_{(r_{i}^{m},\mathbf{z}^{i,m})}\rightarrow \sum\nolimits_{i\geq 1}\delta
_{(r_{i},\mathbf{u}^{i})}$ a.s. \ \rule{0.5em}{0.5em}

\subsubsection{A.s. convergence of versions of $%
(m-F_{1}^{(m)},(F_{2}^{(m)},...))$ to a version of $(\protect\sigma _{I},FI)$%
}

In the following, we keep the notations $r_{i}^{m},$ $\mathbf{z}^{i,m},$ $%
\mathbf{u}^{i,m}$ introduced in the proof above and we recall that we may
suppose that $\sum\nolimits_{i\geq 1}\delta _{(r_{i}^{m},\mathbf{u}^{i,m})}$
and $\sum\nolimits_{i\geq 1}\delta _{(r_{i}^{m},\mathbf{z}^{i,m})}$ converge
to $\sum\nolimits_{i\geq 1}\delta _{(r_{i},\mathbf{u}^{i})}$ a.s. Consider
then some i.i.d. family of $\nu $-homogeneous fragmentations issued from $%
(1,0,...$), say $F^{\text{hom,}(i,j)}$, $i,j\geq 1$, and for each pair $%
(i,j) $, construct from $F^{\text{hom,}(i,j)}$ some $(\tau ,\nu )$%
-fragmentations $F^{(u_{j}^{i,m})}$, $m\geq 1,$ and $F^{(u_{j}^{i})}$,
starting respectively from $u_{j}^{i,m},$ $m\geq 1,$ and $u_{j}^{i}$. Extend
the definition of these processes to $t\in \mathbb{R}^{\ast -}$ by setting $%
F^{(u_{j}^{i,m})}(t)=F^{(u_{j}^{i})}(t):=\mathbf{0}$. Then for $t\geq 0$, let%
\begin{equation}
F^{(i,j),m}(t):=F^{(u_{j}^{i,m})}(t-r_{i}^{m})  \label{9}
\end{equation}%
and set%
\begin{equation}
\overline{\Lambda }^{(m)}(t):=m\prod\nolimits_{r_{i}^{m}\leq t}\left(
1-m^{-1}\sum\nolimits_{j\geq 1}z_{j}^{i,m}\right) .  \label{6}
\end{equation}%
The process $\overline{\Lambda }^{(m)}$ is distributed as $\Lambda ^{(m)}$
since $\sum\nolimits_{j\geq 1}s_{j}^{i}=1$ $\nu $-a.e. for all $i\geq 1$.
The point is then that the process $\overline{F}^{(m)}$ obtained by
considering for each $t\geq 0$ the decreasing rearrangement of the terms $%
\overline{\Lambda }^{(m)}(t)$, $F_{k}^{(i,j),m}(t)$, $i,j,k\geq 1$, is
distributed as $F^{(m)}.$ Furthermore, if $t_{0}<\infty $ is fixed, then
a.s. for $m$ large enough and all $0\leq t\leq t_{0}$, $\overline{F}%
_{1}^{(m)}(t)=\overline{\Lambda }^{(m)}(t)$ and

\begin{equation}
L^{(m)}(t):=(\overline{F}_{2}^{(m)}(t),\overline{F}_{3}^{(m)}(t),...)=%
\{F_{k}^{(i,j),m}(t),i,j,k\geq 1\}^{\downarrow },  \label{11}
\end{equation}%
as noticed in the heuristic description.

Define similarly
\begin{equation*}
F^{(i,j),I}(t):=F^{(u_{j}^{i})}(t-r_{i}),\text{ }t\geq 0\text{.}
\end{equation*}%
The process of the decreasing rearrangements of terms of $F^{(i,j),I}(t)$, $%
i,j\geq 1$, which we still denote by $FI$, is a $(\tau ,\nu ,I)$%
-fragmentation with immigration starting from $\mathbf{0.}$ Also, we still
call $\sigma _{I}(t):=\sum\nolimits_{r_{i}\leq t,j\geq 1}u_{j}^{i},$ $t\geq
0 $. Theorem \ref{Thprincipal} is thus a direct consequence of the following
convergence:

\begin{lemma}
\label{LemmaJoint}$(m-\overline{\Lambda }^{(m)},L^{(m)})\overset{\mathrm{a.s.%
}}{\rightarrow }(\sigma _{I},FI)$ as $m\rightarrow \infty .$
\end{lemma}

To prove this convergence, we shall prove that the following assertions hold
whenever $m_{n}\rightarrow \infty $, and $t_{n}\rightarrow t$, $t_{n}\geq 0$%
, (from now on, we omit the ``a.s.''):

$(A_{a})$ if $t$ is not a jump time of $(\sigma _{I},FI)$, then $(m_{n}-%
\overline{\Lambda }^{(m_{n})}(t_{n}),L^{(m_{n})}(t_{n}))\rightarrow (\sigma
_{I}(t),FI(t))$

$(A_{b})$ if $t$ is a jump time of $(\sigma _{I},FI)$, there exist two
increasing integer-valued sequences (one of them may be finite) $%
\boldsymbol{\varphi }$, $\boldsymbol{\psi }$ such that $\mathbb{N=\{}\varphi
_{n},\psi _{n},n\geq 1\mathbb{\}}$ and

\qquad $(i)$ if $\boldsymbol{\varphi }$ is infinite, then for all sequences $%
s_{\varphi _{n}}\rightarrow t$ s.t. $s_{\varphi _{n}}\geq t_{\varphi _{n}}$,

$\qquad \qquad \qquad (m_{\varphi _{n}}-\overline{\Lambda }^{(m_{\varphi
_{n}})}(s_{\varphi _{n}}),L^{(m_{\varphi _{n}})}(s_{\varphi
_{n}}))\rightarrow (\sigma _{I}(t),FI(t))$

\qquad $(ii)$ if $\boldsymbol{\psi }$ is infinite, then for all sequences $%
s_{\psi _{n}}\rightarrow t$ s.t. $s_{\psi _{n}}\leq t_{\psi _{n}}$,

$\qquad \qquad \qquad (m_{\psi _{n}}-\overline{\Lambda }^{(m_{\psi
_{n}})}(s_{\psi _{n}}),L^{(m_{\psi _{n}})}(s_{\psi _{n}}))\rightarrow
(\sigma _{I}(t-),FI(t-)).$

\noindent According to Lemma \ref{LemmaEthr}, this is sufficient to conclude
that $(m-\overline{\Lambda }^{(m)},L^{(m)})\rightarrow (\sigma _{I},FI)$
with respect to the Skorohod topology. In order to prove these assertions,
we first show two preliminary lemmas.

\begin{lemma}
\textbf{\ }\label{LemmaSko} \textit{Consider a sequence }$a_{n}\rightarrow a$%
\textit{, }$a_{n}\geq 0$\textit{, and let }$F^{\mathrm{\hom }}$\textit{\ be
a }$\nu $-\textit{homogeneous fragmentation starting from }$(1,0,...)$%
\textit{. Let }$F^{(a_{n})}$\textit{\ and }$F^{(a)}$\textit{\ be some }$%
(\tau ,\nu )$\textit{-fragmentations constructed from }$F^{\mathrm{\hom }}$%
\textit{\ starting respectively from }$a_{n},$\textit{\ }$n\geq 0$\textit{,
and }$a$\textit{, and extend these processes to }$t\in \mathbb{R}^{\ast -}$%
\textit{\ by setting }$F^{(a_{n})}(t)=F^{(a)}(t)=\mathbf{0}$. \textit{Then,
whenever }$v_{n}\rightarrow v,$ $v_{n},v\in \mathbb{R},$\textit{\ one has,}

$(a)$\textit{\ if }$v$\textit{\ is not a jump time of }$F^{(a)}$\textit{, }$%
F^{(a_{n})}(v_{n})\rightarrow F^{(a)}(v)$

$(b)$\textit{\ if }$v$\textit{\ is a jump time of }$F^{(a)}$\textit{, there
exist two increasing sequences }$\boldsymbol{\varphi }$\textit{,}$%
\boldsymbol{\psi }$ \textit{such that }$\mathbb{N}=\{\varphi _{n},\psi
_{n},n\geq 1\}$\textit{\ and}

\textit{\ \ \ }$(i)$\textit{\ if }$\boldsymbol{\varphi }$\textit{\ is
infinite and if }$w_{\varphi _{n}}\rightarrow v,$\textit{\ }$w_{\varphi
_{n}}\geq v_{\varphi _{n}},$\textit{\ then }$F^{(a_{\varphi
_{n}})}(w_{\varphi _{n}})\rightarrow F^{(a)}(v)$

$\ \ \ (ii)$\textit{\ if }$\boldsymbol{\psi }$\textit{\ is infinite and if }$%
w_{\psi _{n}}\rightarrow v,$\textit{\ }$w_{\psi _{n}}\leq v_{\psi _{n}},$%
\textit{\ then }$F^{(a_{\psi _{n}})}(w_{\psi _{n}})\rightarrow \nolinebreak
F^{(a)}(v-)$\textit{.}

In particular, when $v=0$, $\boldsymbol{\varphi }$ is the increasing
rearrangement of $\{k:v_{k}\geq 0\}$ and $\boldsymbol{\psi }$ is that of $%
\{k:v_{k}<0\}$.
\end{lemma}

This implies that $F^{(a_{n})}\rightarrow F^{(a)}$\ a.s. with respect to the
Skorohod topology.

\bigskip

\noindent \textbf{Proof.} All the convergences stated in this proof are a.s.
Note that the statement is obvious when $a=0$, since $\sup_{v}\sum_{k\geq
1}F_{k}^{(a_{n})}(v)\leq a_{n}\rightarrow 0.$ Also, when $v_{n}\rightarrow
v<0$, $F^{(a_{n})}(v_{n})=F^{(a)}(v)=0$ for large $n$ and the statement
holds. So, we suppose in the following that $a,a_{n}>0$,\textbf{\ }and $%
v\geq 0$.

To start with, we point out two convergence results when $w_{n}\rightarrow v$%
, $w_{n}\geq 0$. The notations $I^{\hom }$, $T_{x}^{m}$, $m\geq 0$, $x\in
\left( 0,1\right) $, were introduced in Section 1.1.1. First, we claim that $%
T_{x}^{a_{n}}(w_{n})\rightarrow T_{x}^{a}(v)$, provided $\left| I_{x}^{\text{%
hom}}(r)\right| >0$ for all $r\geq 0$ (which occurs for a.e. $x$ $\in \left(
0,1\right) $, since $\nu (\sum\nolimits_{j\geq 1}s_{j}<1)=0$). Indeed,
consider such $x$. If there were a subsequence $(k_{n})_{n\geq 0}$ s.t. $%
\lim_{n\rightarrow \infty }T_{x}^{a_{k_{n}}}(w_{k_{n}})>T_{x}^{a}(v)$ with $%
T_{x}^{a}(v)<\infty $ (the limit may be infinite), then one would have
\begin{equation*}
w_{k_{n}}\geq \int_{0}^{T_{x}^{a_{k_{n}}}(w_{k_{n}})}\mathrm{d}r/\tau
(a_{k_{n}}\left| I_{x}^{\text{hom}}(r)\right|
)>\int_{0}^{T_{x}^{a}(v)+\varepsilon }\mathrm{d}r/\tau (a_{k_{n}}\left|
I_{x}^{\text{hom}}(r)\right| )
\end{equation*}%
for some $\varepsilon >0$ and all $n$ large enough. The latter integral
would then converge to \linebreak $\int_{0}^{T_{x}^{a}(v)+\varepsilon }%
\mathrm{d}r/\tau (a\left| I_{x}^{\text{hom}}(r)\right| )>v$, by dominated
convergence (note that under our assumptions, for $n_{0}$ large enough, the
set $\{a_{k_{n}}\left| I_{x}^{\text{hom}}(r)\right| ,r\leq
T_{x}^{a}(v)+\varepsilon ,n\geq n_{0}\}$ belongs to some compact of $%
(0,\infty )$).\ This would lead to $\lim \inf_{n\rightarrow \infty
}w_{k_{n}}>v$, which is impossible. Similarly, it is not possible that $%
T_{x}^{a_{k_{n}}}(w_{k_{n}})\rightarrow b<T_{x}^{a}(v).$ Hence%
\begin{equation}
T_{x}^{a_{n}}(w_{n})\rightarrow T_{x}^{a}(v)\text{ for a.e. }x\in \left(
0,1\right) .  \label{13}
\end{equation}%
Next, the total mass $M^{(a)}(v)$ can be written as $a\int_{0}^{1}\mathbf{1}%
_{\{T_{x}^{a}(v)<\infty \}}\mathrm{d}x$ and since it is continuous
(Proposition \ref{Propmassecontinue}), $\int_{0}^{1}\mathbf{1}%
_{\{T_{x}^{a}(v)=\infty \}}\mathrm{d}x=0$.$\ $By combining this with (\ref%
{13}), we get
\begin{equation}
M^{(a_{n})}(w_{n})=a_{n}\int_{0}^{1}\mathbf{1}_{\{T_{x}^{a_{n}}(w_{n})<%
\infty \}}\mathbf{1}_{\{T_{x}^{a}(v)<\infty \}}\mathrm{d}x\rightarrow
M^{(a)}(v).  \label{34}
\end{equation}

We are now ready to prove assertion $(a)$. Suppose that $v$ is not a jump
time of $F^{(a)}$. Then for all $x\in (0,1)$, $s\mapsto \left| I_{x}^{\text{%
hom}}(s)\right| $ is continuous at $T_{x}^{a}(v)$ and since $v$ is
necessarily strictly positive, we may suppose that $v_{n}\geq 0$ and apply $%
\left( \ref{13}\right) $. Therefore, $F_{k}^{(a_{n})}(v_{n})\rightarrow
F_{k}^{(a)}(v)$ for all $k\geq 1.$ On the other hand, $M^{(a_{n})}(v_{n})%
\rightarrow M^{(a)}(v)$ by (\ref{34})$.$ Hence $F^{(a_{n})}(v_{n})%
\rightarrow F^{(a)}(v).$

We now turn to $(b)$ and first suppose that $v>0$. That $v$ is a jump time
of $F^{(a)}$\ means that there exists a \textit{unique} interval component
of its interval representation that splits at time $v$.\ More precisely, it
means that there exists a unique interval component, say $I_{v,a}^{\text{hom}%
}$, of the interval representation of $F^{\text{hom}}$ that splits at some
time $T^{a}(v)$ such that $T^{a}(v)=T_{x}^{a}(v)$ for all $x\in I_{v,a}^{%
\text{hom}}$. Moreover, for all $s\leq T^{a}(v)$ and all $x,y\in I_{v,a}^{%
\text{hom}}$, $I_{x}^{\text{hom}}(s)=I_{y}^{\text{hom}}(s)$, which implies
that for all $b,u>0$ and $x\in I_{v,a}^{\text{hom}}$, $T_{x}^{b}(u)\leq
T^{a}(v)\Rightarrow T_{y}^{b}(u)=T_{x}^{b}(u)$ $\forall y\in I_{v,a}^{\text{%
hom}}$.\ This allows us to introduce the increasing sequence $%
\boldsymbol{\psi },$ \textit{independent of }$x\in I_{v,a}^{\text{hom}},$ of
all integers $k$ s.t. $T_{x}^{a_{k}}(v_{k})<T_{x}^{a}(v)$ for some (hence
all) $x\in I_{v,a}^{\text{hom}}$. The increasing\ sequence $%
\boldsymbol{\varphi }$ is then that of all integers $k$ s.t. $%
T_{x}^{a_{k}}(v_{k})\geq T_{x}^{a}(v)$ for some (all) $x\in I_{v,a}$.

Suppose then that $\boldsymbol{\varphi }$ is infinite and let $w_{\varphi
_{n}}\rightarrow v$, $w_{\varphi _{n}}\geq v_{\varphi _{n}}$. On the one
hand, $T_{x}^{a_{\varphi _{n}}}(w_{\varphi _{n}})\geq T_{x}^{a}(v)$, $n\geq
1 $, and the functions $s\mapsto \left| I_{x}(s)\right| $ are
right-continuous when $x\in I_{v,a}^{\text{hom}}.$ On the other hand, the
functions $s\mapsto \left| I_{x}(s)\right| $ are continuous when $x\notin
I_{v,a}^{\text{hom}}$. Therefore, the convergences (\ref{13}) imply that $%
F_{k}^{(a_{\varphi _{n}})}(w_{\varphi _{n}})$ converges to $F_{k}^{(a)}(v),$
$\forall k\geq 1.$ Moreover, $M^{a_{\varphi _{n}}}(w_{\varphi _{n}})$
converges to $M^{a}(v)$ by (\ref{34}) and then, $F^{(a_{\varphi
_{n}})}(w_{\varphi _{n}})$ converges to $F^{(a)}(v).$ Hence $(b)(i).$

Suppose next that $\boldsymbol{\psi }$ is infinite and let $w_{\psi
_{n}}\rightarrow v$, $w_{\psi _{n}}\leq v_{\psi _{n}}$. One has $%
T_{x}^{a_{\psi _{n}}}(w_{\psi _{n}})<T_{x}^{a}(v)$ for all $x\in I_{v,a}^{%
\text{hom}}$, $n\geq 1$. By (\ref{13}), this implies that $F_{k}^{(a_{\psi
_{n}})}(w_{\psi _{n}})$ converges to $F_{k}^{(a)}(v-)$, $\forall k\geq 1$.
Then, using (\ref{34}), we get that $F^{(a_{\psi _{n}})}(w_{\psi _{n}})$
converges to $F^{(a)}(v-)$. Hence $(b)(ii).$

Last, it remains to prove $(b)$ when $v=0.$\ Let here $\boldsymbol{\varphi }$
be the increasing rearrangement of $\{k:v_{k}\geq 0\}$ and $%
\boldsymbol{\psi
}$ the increasing rearrangement of $\{k:v_{k}<0\}.$ If $\boldsymbol{\varphi }
$ is infinite, let $w_{\varphi _{n}}\rightarrow 0$, $w_{\varphi _{n}}\geq
v_{\varphi _{n}}.$\ Then $w_{\varphi _{n}}\geq 0$, and so, according to (\ref%
{13}) and (\ref{34}), $F^{(a_{\varphi _{n}})}(w_{\varphi _{n}})$ converges
to $F^{(a)}(0)$. If $\boldsymbol{\psi }$ is infinite, let $w_{\psi
_{n}}\rightarrow 0$, $w_{\psi _{n}}\leq v_{\psi _{n}}$. Then $F^{(a_{\psi
_{n}})}(w_{\psi _{n}})=0=F^{(a)}(0-).$\ \ \rule{0.5em}{0.5em}

\begin{lemma}
\label{LemmaJoint2}Let $m_{n}\rightarrow \infty $ and $t_{n}\rightarrow t$, $%
t_{n}\geq 0$.

$\mathrm{(i)}$ If $t\notin \{r_{i},i\geq 1\}$, then $m_{n}-\overline{\Lambda
}^{(m_{n})}(t_{n})\rightarrow \sigma _{I}(t)$.

$\mathrm{(ii)}$ If $t=r_{i_{0}}$ for some $i_{0}$, then $m_{n}-\overline{%
\Lambda }^{(m_{n})}(t_{n})$ converges to $\sigma _{I}(t)$ when $t_{n}\geq
r_{i_{0}}^{m_{n}}$ for large $n$'s and it converges to $\sigma _{I}(t-)$
when $t_{n}<r_{i_{0}}^{m_{n}}$ for large $n$'s.
\end{lemma}

\noindent \textbf{Proof. }Recall that $\sum\nolimits_{i}\delta
_{(r_{i}^{m_{n}},\mathbf{z}^{i,m_{n}})}\rightarrow \sum\nolimits_{i}\delta
_{(r_{i},\mathbf{u}^{i})}$ and set $Z^{i,m_{n}}:=\sum\nolimits_{j\geq
1}z_{j}^{i,m_{n}}$, $U^{i}:=\sum\nolimits_{j\geq 1}u_{j}^{i}$.

\textbf{(i) }Take $t^{\prime }>t$ s.t. $t^{\prime }\notin \{r_{i},i\geq 1\}$
and fix $0<\eta <1/2$. One has $\sum\nolimits_{i>k,r_{i}\leq t^{\prime
}}U^{i}<\eta $ for some $k$ large enough and then $\sum%
\nolimits_{i>k,r_{i}^{m_{n}}\leq t^{\prime }}Z^{i,m_{n}}<\eta $ for all $n$
large enough. In particular, all components of these sums are then smaller
than $\eta $. Taking $n$ larger if necessary (so that $m_{n}\geq 1$) one
gets that for all $i>k$, $m_{n}^{-1}Z^{i,m_{n}}\mathbf{1}_{\{r_{i}^{m_{n}}%
\leq t^{\prime }\}}<\eta <1/2$, which implies (using $\left| \ln
(1-x)\right| \leq 2x$ for $0<x\leq 1/2$) that
\begin{equation*}
\left| m_{n}\ln \left( 1-m_{n}^{-1}Z^{i,m_{n}}\right) \mathbf{1}%
_{\{r_{i}^{m_{n}}\leq t^{\prime }\}}\right| \leq 2Z^{i,m_{n}}\mathbf{1}%
_{\{r_{i}^{m_{n}}\leq t^{\prime }\}}\leq 2\eta .
\end{equation*}%
\ When moreover $t_{n}\leq t^{\prime },$%
\begin{equation}
\begin{array}{l}
\sum\nolimits_{i\geq 1}(U^{i}\mathbf{1}_{\left\{ r_{i}\leq t\right\}
}+m_{n}\ln \left( 1-m_{n}^{-1}Z^{i,m_{n}}\right) \mathbf{1}%
_{\{r_{i}^{m_{n}}\leq t_{n}\}}) \\
\leq \sum\nolimits_{i\leq k}(U^{i}\mathbf{1}_{\left\{ r_{i}\leq t\right\}
}+m_{n}\ln \left( 1-m_{n}^{-1}Z^{i,m_{n}}\right) \mathbf{1}%
_{\{r_{i}^{m_{n}}\leq t_{n}\}})+\sum\nolimits_{i>k}(U^{i}\mathbf{1}_{\left\{
r_{i}\leq t^{\prime }\right\} }+2Z^{i,m_{n}}\mathbf{1}_{\{r_{i}^{m_{n}}\leq
t^{\prime }\}}) \\
\leq \sum\nolimits_{i\leq k}(U^{i}\mathbf{1}_{\left\{ r_{i}\leq t\right\}
}+m_{n}\ln \left( 1-m_{n}^{-1}Z^{i,m_{n}}\right) \mathbf{1}%
_{\{r_{i}^{m_{n}}\leq t_{n}\}})+3\eta .%
\end{array}
\label{19}
\end{equation}%
On the other hand, since $t\notin \{r_{i},i\geq 1\},$ $Z^{i,m_{n}}\mathbf{1}%
_{\{r_{i}^{m_{n}}\leq t_{n}\}}\rightarrow U^{i}\mathbf{1}_{\left\{ r_{i}\leq
t\right\} }$ for all $i\geq 1$, or equivalently,%
\begin{equation*}
-m_{n}\ln \left( 1-m_{n}^{-1}Z^{i,m_{n}}\right) \mathbf{1}%
_{\{r_{i}^{m_{n}}\leq t_{n}\}}\rightarrow U^{i}\mathbf{1}_{\left\{ r_{i}\leq
t\right\} }.
\end{equation*}%
Hence the upper bound of $(\ref{19})$ is bounded from above by $4\eta $ for $%
n$ large enough. Therefore $-m_{n}\sum\nolimits_{r_{i}^{m_{n}}\leq t_{n}}\ln
(1-m_{n}^{-1}Z^{i,m_{n}})\ $converges to $\sigma _{I}(t)$ ($%
=\sum\nolimits_{r_{i}\leq t}U^{i}$), which implies that
\begin{equation*}
m_{n}\left( 1-\prod\nolimits_{r_{i}^{m_{n}}\leq t_{n}}\left(
1-m_{n}^{-1}Z^{i,m_{n}}\right) \right) \rightarrow \sigma _{I}(t).
\end{equation*}

\textbf{(ii) }If $t=r_{i_{0}}$ and $t_{n}\geq r_{i_{0}}^{m_{n}}$ for $n$
large enough, then, for all $i\geq 1,$ $Z^{i,m_{n}}\mathbf{1}%
_{\{r_{i}^{m_{n}}\leq t_{n}\}}$ converges to $U^{i}\mathbf{1}_{\{r_{i}\leq
t\}}$ and one concludes exactly as above. Now, if $t_{n}<r_{i_{0}}^{m_{n}}$
for large $n$'s, $Z^{i_{0},m_{n}}\mathbf{1}_{\{r_{i_{0}}^{m_{n}}\leq
t_{n}\}} $ converges to $U^{i_{0}}\mathbf{1}_{\{r_{i_{0}}<t\}}$ and still, $%
Z^{i,m_{n}}\mathbf{1}_{\{r_{i}^{m_{n}}\leq t_{n}\}}$ converges to $U^{i}%
\mathbf{1}_{\{r_{i}\leq t\}}$ for $i\neq i_{0}$. The conclusion then follows
by replacing $\mathbf{1}_{\{r_{i_{0}}\leq t\}}$\ by $\mathbf{1}%
_{\{r_{i_{0}}<t\}}$ in the proof above. \ \rule{0.5em}{0.5em}

\bigskip

Now let $m_{n}\rightarrow \infty $ and $t_{n}\rightarrow t$, $t_{n}\geq 0$.
We are ready to prove assertions $(A_{a})$ and $(A_{b})$.\bigskip

\noindent \textbf{Proof of Lemma \ref{LemmaJoint}.} The proof is split into
two parts, according to whether $t$ is, or not, a jump time of $FI$. It
strongly relies on the convergence $\sum\nolimits_{i\geq 1}\delta
_{(r_{i}^{m_{n}},\mathbf{u}^{i,m_{n}})}\rightarrow \sum\nolimits_{i\geq
1}\delta _{(r_{i},\mathbf{u}^{i})}$.

\textbf{1}. If $t$ is not a jump time of $FI$, then $t-r_{i}$ is not a jump
time of $F^{(u_{j}^{i})}$, $\forall i,j\geq 1$ (in particular, $t\notin
\{r_{i},i\geq 1\}$). When $t-r_{i}>0$, applying Lemma \ref{LemmaSko} $(a)$
to the sequences $a_{n}=u_{j}^{i,m_{n}}$, $a=u_{j}^{i}$, $%
v_{n}=t_{n}-r_{i}^{m_{n}}$ and $v=t-r_{i}$, one gets that $%
F^{(i,j),m_{n}}(t_{n})\rightarrow F^{(i,j),I}(t)$, $\forall j\geq 1$.
Clearly, such convergence also holds when $t<r_{i}$, since then $%
F^{(i,j),m_{n}}(t_{n})=\mathbf{0}$ for $n$ large enough. Fix next some $\eta
,\varepsilon >0$ and then some $k$ s.t. $\sum\nolimits_{i+j>k}u_{j}^{i}%
\mathbf{1}_{\{r_{i}\leq t+\varepsilon \}}<\eta $. For $n$ large enough, $%
\sum\nolimits_{i+j>k}u_{j}^{i,m_{n}}\mathbf{1}_{\{r_{i}^{m_{n}}\leq
t_{n}\}}<\eta $, and therefore%
\begin{equation*}
\sum\nolimits_{i,j\geq 1}d(F^{(i,j),m_{n}}(t_{n}),F^{(i,j),I}(t))\leq
\sum\nolimits_{i+j\leq k}d(F^{(i,j),m_{n}}(t_{n}),F^{(i,j),I}(t))+2\eta .
\end{equation*}%
So, the right hand side of this inequality is smaller than $3\eta $ for $n$
large enough, i.e. \linebreak $\sum\nolimits_{i,j\geq
1}d(F^{(i,j),m_{n}}(t_{n}),F^{(i,j),I}(t))\rightarrow 0$. Then, by Lemma \ref%
{Lemmadecresase}, one concludes that $L^{(m_{n})}(t_{n})\rightarrow FI(t)$.\
On the other hand, Lemma \ref{LemmaJoint2} (i) implies that $m_{n}-\overline{%
\Lambda }^{(m_{n})}(t_{n})\rightarrow \sigma _{I}(t).$ Hence we have
assertion $(A_{a})$.

\textbf{2}. Now assume that $t$ is a jump time of $FI$. Our goal is to
construct some increasing sequences $\boldsymbol{\varphi }$ and $%
\boldsymbol{\psi }$, $\mathbb{N}=\{\varphi _{n},\psi _{n},n\geq 1\},$ such
that assertions $(A_{b})(i)$ and $(A_{b})(ii)$ hold. For all $i,j\geq 1$,
let $\mathcal{J}^{(i,j)}$ denote the set of \textit{strictly positive} jump
times of $F^{(u_{j}^{i})}$. Since the process of total mass of this
fragmentation is continuous (Proposition \ref{Propmassecontinue}), it only
jumps when a fragment splits. And then, since the $F^{\text{hom,}(i,j)}$ are
constructed from independent Poisson point processes, independent of $%
((r_{i},\mathbf{u}^{i}),i\geq 1)$, (a.s.) the $\mathcal{J}^{(i,j)}$'s are
pairwise disjoint and disjoint from $\{r_{i},i\geq 1\}$. Also, every $%
F^{(u_{j}^{i})}$ jumps at $0$ (we recall that these processes are defined on
$\mathbb{R}$) and therefore the set of jump times of $F^{(u_{j}^{i})}$ is $%
\mathcal{J}^{(i,j)}\cup \{0\}$, $i,j\geq 1$. So, if $t$ is a jump time of $%
FI $:

$\bullet $ either $t-r_{i_{0}}\in \mathcal{J}^{(i_{0},j_{0})}$ for some
(unique) pair $(i_{0},j_{0}).$ Then, one can apply Lemma \ref{LemmaSko}$(b)$
to $a_{n}=u_{j_{0}}^{i_{0},m_{n}}$, $a=u_{j_{0}}^{i_{0}}$, $%
v_{n}=t_{n}-r_{i_{0}}^{m_{n}}$ and $v=t-r_{i_{0}}$. Let $%
\boldsymbol{\varphi
}$ and $\boldsymbol{\psi }$ be the sequences that appear in this statement%
\textit{\ }and first, suppose that $\boldsymbol{\varphi }$ is infinite.
Consider then some sequence $s_{\varphi _{n}}\rightarrow t$, $s_{\varphi
_{n}}\geq t_{\varphi _{n}}$, and apply Lemma \ref{LemmaSko}$(b)(i)$ to $%
w_{\varphi _{n}}=s_{\varphi _{n}}-r_{i_{0}}^{m_{\varphi _{n}}}$. One obtains
that $F^{(i_{0},j_{0}),m_{\varphi _{n}}}(s_{\varphi _{n}})$ converges to $%
F^{(i_{0},j_{0}),I}(t)$. On the other hand, by Lemma \ref{LemmaSko}$(a)$, $%
F^{(i,j),m_{\varphi _{n}}}(s_{\varphi _{n}})\ $converges to $F^{(i,j),I}(t)$
for all $(i,j)\neq (i_{0},j_{0})$ since $t-r_{i}\notin \mathcal{J}%
^{(i,j)}\cup \{0\}$. Hence, as in \textbf{1}., we get that $%
\sum\nolimits_{i,j\geq 1}d(F^{(i,j),m_{\varphi _{n}}}(s_{\varphi
_{n}}),F^{(i,j),I}(t))$ tends to $0$, and then, by Lemma \ref{Lemmadecresase}%
, that $L^{m_{\varphi _{n}}}(s_{\varphi _{n}})$ converges to $FI(t).$
Moreover, $m_{\varphi _{n}}-\overline{\Lambda }^{(m_{\varphi
_{n}})}(s_{\varphi _{n}})$ converges to $\sigma _{I}(t),$ by Lemma \ref%
{LemmaJoint2} (i) since $t\notin \{r_{i},i\geq 1\}$. Hence $(A_{b})(i).$

Similarly, supposing that $\boldsymbol{\psi }$ is infinite and $s_{\psi
_{n}}\rightarrow t$, $s_{\psi _{n}}\leq t_{\psi _{n}}$, one gets, by
applying Lemma \ref{LemmaSko}$(b)(ii)$, that $L^{m_{\psi _{n}}}(s_{\psi
_{n}})\rightarrow FI(t-)$. Moreover, $m_{\psi _{n}}-\overline{\Lambda }%
^{(m_{\psi _{n}})}(s_{\psi _{n}})$ converges to $\sigma _{I}(t)=\sigma
_{I}(t-),$ still by Lemma \ref{LemmaJoint2} (i) since $t\notin \{r_{i},i\geq
1\}$. Hence $(A_{b})(ii)$ and then $(A_{b})$.

$\bullet $ or $t\in \{r_{i},i\geq 1\},$ say $t=r_{i_{0}}$. For $i\neq i_{0}$%
, $t-r_{i}\notin \mathcal{J}^{(i,j)}\cup \{0\}$ and therefore, as explain
above, $F^{(i,j),m_{n}}(s_{n})$ converges to $F^{(i,j),I}(t)=F^{(i,j),I}(t-)$%
, for all sequences $s_{n}\rightarrow t$. Let then $\boldsymbol{\varphi }$
be the increasing sequence of integers $k$ such that $t_{k}\geq
r_{i_{0}}^{m_{k}}$ and $\boldsymbol{\psi }$ be the increasing sequence of
integers $k$ such that $t_{k}<r_{i_{0}}^{m_{k}}$. When $\boldsymbol{\varphi }
$ is infinite and $s_{\varphi _{n}}\rightarrow t=r_{i_{0}}$, $s_{\varphi
_{n}}\geq t_{\varphi _{n}}$, one has, by Lemma \ref{LemmaSko}(b)(i), that $%
F^{(i_{0},j),m_{\varphi _{n}}}(s_{\varphi _{n}})$ converges to $%
F^{(i_{0},j),I}(t)$, $\forall j\geq 1$. Together with the fact that $%
F^{(i,j),m_{\varphi _{n}}}(s_{\varphi _{n}})$ converges to $F^{(i,j),I}(t)$
for $i\neq i_{0}$, $j\geq 1$, we obtain, as in \textbf{1.}, that $%
L^{(m_{\varphi _{n}})}(s_{\varphi _{n}})$ converges to $FI(t)$. On the other
hand, $m_{\varphi _{n}}-\overline{\Lambda }^{(m_{\varphi _{n}})}(s_{\varphi
_{n}})$ converges to $\sigma _{I}(t),$ by Lemma \ref{LemmaJoint2} (ii).
Hence assertion $(A_{b})(i)$. Now, if $\boldsymbol{\psi }$ is infinite, let $%
s_{\psi _{n}}\rightarrow t$, $s_{\psi _{n}}\leq t_{\psi _{n}}$. Clearly, $%
F^{(i_{0},j),m_{\psi _{n}}}(s_{\psi _{n}})=0=F^{(i_{0},j),I}(t-)$, $\forall
j,n\geq 1$. Moreover $F^{(i,j),m_{\psi _{n}}}(s_{\psi _{n}})$ converges to $%
F^{(i,j),I}(t-)$ for $i\neq i_{0}$, $j\geq 1$, and therefore $L^{(m_{\psi
_{n}})}(s_{\psi _{n}})$ tends to $FI(t-)$. At last, $m_{\psi _{n}}-\overline{%
\Lambda }^{(m_{\psi _{n}})}(s_{\psi _{n}})$ converges to $\sigma _{I}(t-)$
by Lemma \ref{LemmaJoint2} (ii). Hence assertion $(A_{b})(ii)$. \ \rule%
{0.5em}{0.5em}

\textbf{\newline
Remark. }The convergence in law of $m-\Lambda ^{(m)}$ (and a fortiori of $%
m-F_{1}^{(m)}$) to some $\gamma $-stable subordinator $\sigma $, $\gamma \in
(0,1)$, actually holds as soon as $\varphi _{\nu }$ varies regularly at $%
\infty $ with index $-\gamma $ and $\tau (m)\sim C\varphi _{\nu }(m)$, $C>0$%
. Very roughly, the point is either to check that the regular variation
assumptions imply that the measures $\sum\nolimits_{i\geq 1}\delta
_{(r_{i}^{m},\sum\nolimits_{j\geq 1}u_{j}^{i,m})}$ converge in distribution
to some Poisson point measure $\sum\nolimits_{i\geq 1}\delta
_{(r_{i},x^{i})} $ on $\left[ 0,\infty \right) \times \mathbb{R}^{+},$ where
$((r_{i},x^{i}),i\geq 1)$ is a PPP with intensity $C^{\prime }x^{-\gamma -1}$%
d$x$, $x>0$. This will lead to some result identical to Lemma \ref%
{LemmaJoint2} (replacing there $\sigma _{I}$ by $\sigma $). Or to use
classical results on convergence of subordinators (see e.g.\cite%
{JacodShiryaev}) and, again, regular variation theorems.

\subsection{Proof of Theorem \ref{Proposition_infinie}}

We still use the notations $\Lambda ^{(m)},$ $((t_{i}^{m},\mathbf{s}%
^{i}),i\geq 1)$ and $((r_{i},\mathbf{u}^{i},i\geq 1)$ introduced in the
previous subsection, and we suppose that $\tau $\ varies regularly at $%
\infty ,$\ and that $\varphi _{\nu }(m)\nu _{m}\rightarrow I$\ as $%
m\rightarrow \infty $.

\bigskip

\noindent \textbf{Lemma \ref{Lemmamesures} bis} $\sum\nolimits_{i\geq
1}\delta _{(t_{i}^{m}(\tau /\varphi _{\nu })(m),(\Lambda
^{(m)}(t_{i}^{m}-)s_{j}^{i})_{j\geq 2})}\overset{\mathrm{law}}{\rightarrow }%
\sum\nolimits_{i\geq 1}\delta _{(r_{i},\mathbf{u}^{i})}.$

\bigskip

The proof of this result is similar to that of Lemma \ref{Lemmamesures}, so
we only give the main lines.\bigskip

\noindent \textbf{Proof. }Here we replace the measure $\mu _{m}$\ introduced
in the proof of Lemma \ref{Lemmamesures} by $\widetilde{\mu }%
_{m}:=\sum\nolimits_{i\geq 1}\delta _{((t_{i}/\varphi _{\nu
}(m),(ms_{j}^{i})_{j\geq 2})},$ which, under the assumption $\varphi _{\nu
}(m)\nu _{m}\rightarrow I,$ converges in distribution to $%
\sum\nolimits_{i\geq 1}\delta _{(r_{i},\mathbf{u}^{i})}$. Consider versions
of these measures such that the a.s. convergence holds. Then let $\widetilde{%
\sigma }^{m}$ be a permutation such that $t_{\widetilde{\sigma }%
^{m}(i)}/\varphi _{\nu }(m)\rightarrow r_{i}$ and $(ms_{j}^{\widetilde{%
\sigma }^{m}(i)})_{j\geq 2}\rightarrow \mathbf{u}^{i},\ $and define from
this permutation, $\widetilde{r}_{i}^{m},$ $\widetilde{\mathbf{u}}^{i,m}$
and $\widetilde{\mathbf{z}}^{i,m},$ exactly as $r_{i}^{m},$ $\mathbf{u}%
^{i,m} $ and $\mathbf{z}^{i,m}$ were defined from $\sigma ^{m}$ by formula $%
\left( \ref{3}\right) $. As in the proof of Lemma \ref{Lemmamesures}, note
that $\widetilde{\mathbf{u}}^{i,m}\rightarrow \mathbf{u}^{i},$ since $t_{%
\widetilde{\sigma }^{m}(i)}\rightarrow 0$ since $\varphi _{\nu
}(m)\rightarrow 0$ (this convergence to $0$ is due to the fact that $\varphi
_{\nu }$ varies regularly with some index in $(-1,0)$, since $\varphi _{\nu
}(m)\nu _{m}$ converges, see Lemma \ref{Lemma 6}). Also, since $\tau $\
varies regularly at $\infty $, the Potter's bounds (\ref{10}) hold and then
one checks that $\widetilde{r}_{i}^{m}(\tau /\varphi _{\nu })(m)\rightarrow
r_{i}$. The rest of the proof is now very similar to that of Lemma \ref%
{Lemmamesures}. \ \rule{0.5em}{0.5em}

\bigskip

One may suppose that $\sum\nolimits_{i\geq 1}\delta _{(\widetilde{r}%
_{i}^{m}(\tau /\varphi _{\nu })(m),\widetilde{\mathbf{u}}^{i,m})}$ and $%
\sum\nolimits_{i\geq 1}\delta _{(\widetilde{r}_{i}^{m}(\tau /\varphi _{\nu
})(m),\widetilde{\mathbf{z}}^{i,m})}$ converge to\linebreak\ $%
\sum\nolimits_{i\geq 1}\delta _{(r_{i},\mathbf{u}^{i})}$ a.s. Let then$\,%
\widetilde{F}^{(i,j),m}$, $\widetilde{\Lambda }^{(m)}$ and $\widetilde{L}%
^{(m)}$ be defined from $\widetilde{r}_{i}^{m},\widetilde{\mathbf{u}}^{i,m},%
\widetilde{\mathbf{z}}^{i,m},$ $i\geq 1,$ $m\geq 0,$ by formulas similar to (%
\ref{9}), (\ref{6}) and $\left( \ref{11}\right) $. Also, let $\widetilde{F}%
^{(m)}$ obtained by considering for each $t\geq 0$ the decreasing
rearrangement of the terms $\widetilde{\Lambda }^{(m)}(t)$, $\widetilde{F}%
_{k}^{(i,j),m}(t)$, $i,j,k\geq 1$, and note that $\widetilde{\Lambda }^{(m)}%
\overset{\text{law}}{=}\Lambda ^{(m)}$ and $\widetilde{F}^{(m)}\overset{%
\text{law}}{=}F^{(m)}$. We should point out that contrary to what happens
when $\tau (m)\nu _{m}$ converges to a non-trivial limit, $\widetilde{%
\Lambda }^{(m)}$ and $\widetilde{F}_{1}^{(m)}$ do not necessarily coincide
on $[0,t_{0}]$ for large $m$'s under the assumptions of Theorem \ref%
{Proposition_infinie}. However $\widetilde{\Lambda }^{(m)}((\varphi _{\nu
}/\tau )(m)\cdot )$ and $\widetilde{F}_{1}^{(m)}((\varphi _{\nu }/\tau
)(m)\cdot )$ do coincide on $[0,t_{0}]$ for large $m$'s and that is all we
need for the proof.

Let $m_{n}\rightarrow \infty $ and $t_{n}\rightarrow t$. By imitating the
proof of Lemma \ref{LemmaJoint2} one easily obtains\bigskip

\noindent \textbf{Lemma \ref{LemmaJoint2} bis }$\mathrm{(i)}$ \textit{If }$%
t\notin \{r_{i},i\geq 1\}$\textit{, then }$m_{n}-\widetilde{\Lambda }%
^{(m_{n})}((\varphi _{\nu }/\tau )(m_{n})t_{n})\rightarrow \sigma _{I}(t)$%
\textit{.}

$\mathrm{(ii)}$ \textit{If }$t=r_{i_{0}}$\textit{, then }$m_{n}-\widetilde{%
\Lambda }^{(m_{n})}((\varphi _{\nu }/\tau )(m_{n})t_{n})$\textit{\ converges
to }$\sigma _{I}(t)$\textit{\ when }$(\varphi _{\nu }/\tau )(m_{n})t_{n}\geq
\widetilde{r}_{i_{0}}^{m_{n}}$ \textit{for }$n$\textit{\ large enough and it
converges to }$\sigma _{I}(t-)$\textit{\ when }$(\varphi _{\nu }/\tau
)(m_{n})t_{n}<r_{i_{0}}^{m_{n}}$ \textit{for }$n$\textit{\ large enough.}

\bigskip

\noindent \textbf{Proof of Theorem \ref{Proposition_infinie} (i)}. Let $I(t)$
be the decreasing rearrangement $\{u_{j}^{i},j\geq 1,r_{i}\leq
t\}^{\downarrow }$. Our goal is to show that
\begin{equation}
\left( m-\widetilde{\Lambda }^{(m)}\left( (\varphi _{\nu }/\tau )(m)\cdot
\right) ,\widetilde{L}^{(m)}\left( (\varphi _{\nu }/\tau )(m)\cdot \right)
\right) \overset{\text{a.s.}}{\rightarrow }\left( \sigma _{I},(I(t),t\geq
0)\right)  \label{12}
\end{equation}%
when $(\varphi _{\nu }/\tau )(m)\rightarrow 0$. Under this assumption, $%
(\varphi _{\nu }/\tau )(m_{n})\left( t_{n}-\widetilde{r}_{i}^{m_{n}}(\tau
/\varphi _{\nu })(m_{n})\right) \ $converges to $0$. Hence $\widetilde{F}%
^{(i,j),m_{n}}((\varphi _{\nu }/\tau )(m_{n})t_{n})\ $converges to $%
u_{j}^{i} $ when $(\varphi _{\nu }/\tau )(m_{n})t_{n}\geq \widetilde{r}%
_{i}^{m_{n}}$ for large $n$'s and it reaches $0$ when $(\varphi _{\nu }/\tau
)(m_{n})t_{n}<\widetilde{r}_{i}^{m_{n}}$ for large $n$'s. Recalling Lemma %
\ref{LemmaJoint2} bis, it is then easy to adapt the proof of Lemma \ref%
{LemmaJoint} to obtain the required convergence (\ref{12}). Note that the
only jump times of the limit process are the $r_{i}$'s, which makes the
proof shorter than that of Lemma \ref{LemmaJoint}. \ \rule{0.5em}{0.5em}

\bigskip

To prove Theorem \ref{Proposition_infinie} (ii), we need the following lemma.

\begin{lemma}
\label{LemmaCV2}\textit{Consider a sequence }$a_{n}\rightarrow a$\textit{, }$%
a_{n}\geq 0$\textit{, and let }$F^{\mathrm{\hom }}$\textit{\ be a
homogeneous }$\nu $\textit{-fragmentation starting from }$(1,0,...)$\textit{%
. Let }$F^{(a_{n})}$\textit{\ be some }$(\tau ,\nu )$\textit{-fragmentations
constructed from }$F^{\mathrm{\hom }},$\textit{\ starting respectively from }%
$a_{n},$\textit{\ }$n\geq 0$\textit{, and let }$t_{n}\rightarrow \infty .$%
\textit{\ Then, if }the fragmentation loses mass to dust,\textit{\ }$%
F^{(a_{n})}(t_{n})\rightarrow \mathbf{0}$ a.s.
\end{lemma}

\noindent \textbf{Proof. }As in the proof of Lemma \ref{LemmaSko} we may
suppose that $a>0$. Then, since the fragmentation $F^{(a)}$ loses mass,
every $x$ falls into the dust after a finite time, i.e. $\int_{0}^{\infty
}1/\tau (aI_{x}^{\text{hom}}(r))$d$r<\infty .$ It is then easy, using
dominated convergence and the fact that $\tau $ is monotone near $0$ (hence
necessarily non-increasing here, because of the loss of mass), that $%
\int_{0}^{\infty }1/\tau (a_{n}I_{x}^{\text{hom}}(r))$d$r\leq C<\infty $ for
large $n$'s. Hence $T_{x}^{a_{n}}(t_{n})=\infty $ for large $n$'s, and
therefore, $M^{(a_{n})}(t_{n})=a_{n}\int_{0}^{1}\mathbf{1}%
_{\{T_{x}^{a_{n}}(t_{n})<\infty \}}\mathrm{d}x$ converges to $0$. \ \rule%
{0.5em}{0.5em}

\bigskip

\noindent \textbf{Proof of Theorem \ref{Proposition_infinie} (ii)}. Suppose
that $(\varphi _{\nu }/\tau )(m)\rightarrow \infty $ and fix $t\geq 0$. When
$t>r_{i}$, $(\varphi _{\nu }/\tau )(m)\left( t-(\tau /\varphi _{\nu })(m)%
\widetilde{r}_{i}^{m}\right) \rightarrow \infty $ and then, according to the
previous lemma,\linebreak $\widetilde{F}^{(i,j),m}((\varphi _{\nu }/\tau
)(m)t)\rightarrow \mathbf{0}$. When $t<r_{i}$, $\widetilde{F}%
^{(i,j),m}((\varphi _{\nu }/\tau )(m)t)=\mathbf{0}$ for $m$ large enough.
From this, we deduce that for all $t\notin \{r_{i},i\geq 1\}$, $\widetilde{L}%
^{(m)}\left( (\varphi _{\nu }/\tau )(m)t\right) \rightarrow \mathbf{0}$.
Furthermore, $m-\widetilde{\Lambda }^{(m)}\left( (\varphi _{\nu }/\tau
)(m)t\right) \rightarrow \sigma _{I}(t)$ according to Lemma \ref{LemmaJoint2}
bis. So, if we consider some finite sequence of deterministic times $%
t_{1},...,t_{k}$, we know that (a.s.) these times are not in $\{r_{i},i\geq
1\}$, and therefore that the convergences of $(m-\widetilde{\Lambda }%
^{(m)}\left( (\varphi _{\nu }/\tau )(m)t_{l}\right) ,\widetilde{L}%
^{(m)}\left( (\varphi _{\nu }/\tau )(m)t_{l}\right) )$ to $(\sigma
_{I}(t_{l}),\mathbf{0})$, $1\leq l\leq k$, hold simultaneously. Hence the
convergence in the finite dimensional sense. \ \rule{0.5em}{0.5em}

\bigskip

Let us point out that the convergence of $\widetilde{L}^{(m)}\left( (\varphi
_{\nu }/\tau )(m)\cdot \right) $ to $\mathbf{0}$ in the Skorohod sense does
not hold when $(\varphi _{\nu }/\tau )(m)\rightarrow \infty $. Indeed,
consider some $i$ such that $u_{1}^{i}>0$ (such $i$ exists since $%
I(l_{1}^{\downarrow })\neq 0$) and set $t_{m}:=(\tau /\varphi _{\nu })(m)%
\widetilde{r}_{i}^{m}$, $m\geq 0$. Then $\widetilde{F}^{(i,1),m}((\varphi
_{\nu }/\tau )(m)t_{m})$ converges to $u_{1}^{i}\neq 0$ and consequently $%
\widetilde{L}^{(m)}\left( (\varphi _{\nu }/\tau )(m)t_{m}\right)
\nrightarrow \mathbf{0}$. Therefore, assertion (a) of Lemma \ref{LemmaEthr}
is not satisfied.

\section{Small times asymptotics in the self-similar cases}

We are now looking at the small times asymptotics of $F^{(1)}$ when $\tau
(m)=m^{\alpha }$, $\alpha \in \mathbb{R}$. When $\nu (l_{1,\leq
1}^{\downarrow })<\infty ,$ a particle waits a positive time before
splitting and $F^{(1)}(\varepsilon )\overset{\text{\textrm{a.s.}}}{=}%
(1,0,...)$ for $\varepsilon $ small enough. So the interesting case to study
here is $\nu (l_{1,\leq 1}^{\downarrow })=\infty ,$ which is supposed in the
rest of this section. In that aim, introduce the function
\begin{equation*}
\varphi _{\nu }^{-1}(\varepsilon ):=\inf \{m:\varphi _{\nu }(m)<\varepsilon
\},
\end{equation*}%
which is well defined in a neighborhood of $0$ since $\nu (l_{1,\leq
1}^{\downarrow })=\infty $, and recall that under the assumption $\varphi
_{\nu }(m)\nu _{m}\rightarrow I$, the function $\varphi _{\nu }$ is
regularly varying at $\infty $ (with index $-\gamma _{\nu }$). Classical
results on regular variation (see \cite{bgt}) then implies that $\varphi
_{\nu }^{-1}$ is also regularly varying (at $0$) and $\varphi _{\nu }\circ
\varphi _{\nu }^{-1}(\varepsilon )\sim \varepsilon $ when $\varepsilon
\rightarrow 0$. Thanks to the self-similarity of $F$, one then obtains the
following Corollary \ref{Corosmallbeh} by

\vspace{0.2cm}- applying Theorem \ref{Thprincipal} when $m^{-\alpha }\varphi
_{\nu }(m)\rightarrow \ell \in (0,\infty )$ to the process $F^{((\varepsilon
\ell ^{-1})^{1/\alpha })}$, and then using that a fragmentation with
immigration process with parameters $(\alpha ,\ell \nu ,I)$ is distributed
as $FI^{\ell ^{-1}}(\ell \cdot )$ where $FI^{\ell ^{-1}\text{ }}$denotes a
fragmentation with immigration $(\alpha ,\nu ,\ell ^{-1}I)$

- applying Theorem \ref{Proposition_infinie} when $m^{-\alpha }\varphi _{\nu
}(m)\rightarrow \ell \in \{0,\infty \}$ to the process $F^{(\varphi _{\nu
}^{-1}(\varepsilon ))}.$\vspace{0.2cm}

\noindent By convention, when $\ell =\infty $, a $\left( \alpha ,\ell \nu
,I\right) $ fragmentation with immigration $FI$ is a process constantly
null, $FI(t)=\mathbf{0},$ $\forall t\geq 0$, but the subordinator $\sigma
_{I}$ of total mass of immigrants is still non-trivial and constructed from
the measure $I$. Roughly, this corresponds to the case where particles
immigrate and vanish immediately.

\begin{corollary}
\label{Corosmallbeh}Suppose that $\varphi _{\nu }(m)\nu _{m}\rightarrow I$
and $m^{-\alpha }\varphi _{\nu }(m)\rightarrow \ell \in \left[ 0,\infty %
\right] $ as $m\rightarrow \infty ,$ and let $FI$ be a self-similar
fragmentation with immigration with parameters $\left( \alpha ,\ell \nu
,I\right) ,$ starting from $\mathbf{0.}$ Then,%
\begin{equation*}
\varphi _{\nu }^{-1}(\varepsilon )(1-F_{1}^{(1)}(\varepsilon \cdot
),(F_{2}^{(1)}(\varepsilon \cdot ),F_{3}^{(1)}(\varepsilon \cdot ),...))%
\overset{\mathrm{law}}{\rightarrow }(\sigma _{I},FI)\text{ as }\varepsilon
\rightarrow 0,
\end{equation*}%
where the convergence holds with respect to the Skorohod topology when $\ell
<\infty $ and in the finite-dimensional sense when $\ell =\infty .$
\end{corollary}

Remark that the fragmentation with immigration process that arises in this
limit is $\gamma _{\nu }$-\textit{self-similar} (as a consequence of the $%
\gamma _{\nu }$-self-similarity of $I$ stated in Lemma \ref{Lemma 6}), i.e.%
\begin{equation*}
(FI(at),t\geq 0)\overset{\text{law}}{=}(a^{1/\gamma _{\nu }}FI(t),t\geq 0)%
\text{ for all }a>0.
\end{equation*}

Bertoin \cite{bertafrag02} proves that large times behavior of self-similar
fragmentations differs significantly according as $\alpha <0$, $\alpha =0$
or $\alpha >0$. The above corollary shows that the rules are quite different
for small times behavior: the convergence rate only depends on $\nu $ and
then the form of the limit only depends on the position of $\alpha $ with
respect to $\gamma _{\nu }$. The index $\alpha =-\gamma _{\nu }$ is the only
one for which the limit may be a non-trivial fragmentation with immigration
and this occurs if and only if $\varphi _{\nu }(m)$ behaves as a power
function as $m\rightarrow \infty $. This suggests that this index is in some
sense more natural than the others.

However, the limit is also non-trivial when $\alpha >-\gamma _{\nu }.$ In
particular, Corollary \ref{Corosub}, in this self-similar setting, says that
if $\nu $ is binary and if $\varphi _{\nu }$ varies regularly at $\infty $
with some index $-\gamma _{\nu }\in (-1,0)$, then, as soon as $\alpha
>-\gamma _{\nu }$,
\begin{equation*}
\varphi _{\nu }^{-1}(\varepsilon )\left( 1-F_{1}^{(1)}(\varepsilon \cdot
),F_{2}^{(1)}(\varepsilon \cdot ),F_{3}^{(1)}(\varepsilon \cdot ),...\right)
\overset{\mathrm{law}}{\rightarrow }(\sigma ,\Delta _{1},\Delta _{2},...)%
\text{ as }\varepsilon \rightarrow 0
\end{equation*}%
where $\sigma $ is a stable subordinator with Laplace exponent $\Gamma
(1-\gamma _{\nu })q^{\gamma _{\nu }}$ and $(\Delta _{1}(t),\Delta
_{2}(t),...)$ the decreasing sequence of its jumps before time $t,$ $t\geq 0$%
. This completes a result of Berestycki \cite{BerestyckiThese} who shows
that
\begin{equation*}
\varphi _{\nu }^{-1}(\varepsilon )\left( F_{2}^{(1)}(\varepsilon
),F_{3}^{(1)}(\varepsilon ),...\right) \overset{\mathrm{law}}{\rightarrow }%
(\Delta _{1}(1),\Delta _{2}(1),...)
\end{equation*}%
when $\alpha \geq 0$, $\nu $ is binary and $\varphi _{\nu }$ varies
regularly at $\infty $. He also investigates the behavior of $%
F_{2}^{(1)}(\varepsilon )$ near $0$ for all measures $\nu $ and $\alpha \geq
0,$ and obtains that $F_{2}^{(1)}(\varepsilon )\sim R(\varepsilon )$ a.s.
where $R$ is the record process of a PPP with intensity $\nu (s_{2}\in
\mathrm{d}x).$

We also refer to Miermont and Schweinsberg \cite{MiermontSchweinsberg} for
some specific examples.

\bigskip

\noindent \textbf{Total mass behavior. }In the self-similar setting, the
total mass $M^{(1)}(t)=\sum\nolimits_{i\geq 1}F_{i}^{(1)}(t)$ of macroscopic
particles present at time $t$ is non-constant if and only if $\alpha <0$. A
consequence of Corollary \ref{Corosmallbeh} is that the behavior near $0$ of
the mass $1-M^{(1)}$ is then specified as follows.

\begin{corollary}
Under the assumptions of Corollary \ref{Corosmallbeh}, as $\varepsilon
\rightarrow 0$,
\begin{equation*}
\varphi _{\nu }^{-1}(\varepsilon )(1-M^{(1)}(\varepsilon \cdot ))\overset{%
\mathrm{law}}{\rightarrow }\sigma _{I}-M_{FI}
\end{equation*}%
where $M_{FI}(t)=\sum_{j\geq 1}FI_{j}(t),$ $t\geq 0$, (again, the
convergence holds with respect to the Skorohod topology when $\ell <\infty $
and in the finite-dimensional sense when $\ell =\infty $). In particular,
the limit is equal to $\sigma _{I}$ when $\ell =\infty ,$ is $\mathbf{0}$
when $\ell =0$, and is non-trivial when $0<\ell <\infty .$
\end{corollary}

Note that when $\alpha >-\gamma _{\nu }$, the limit $\ell $ equals $0$ and
so the speed of convergence of $1-M^{(1)}(\varepsilon )$ to $0$ is faster
than $1/\varphi _{\nu }^{-1}(\varepsilon )$. When $-\gamma _{\nu }<\alpha <0$%
, one can obtain a lower bound for this speed by using Theorem 4 of \cite%
{HaasMiermont}, which implies that for all $\gamma <-\alpha $, there exists
a positive constant $C_{\gamma }$ such that $1-M^{(1)}(\varepsilon )\geq
C_{\gamma }\varepsilon ^{1/\gamma }$, $\forall \varepsilon >0$.

\section{Underlying continuum random trees}

In this section $\tau (m)=m^{\alpha }$ with $\alpha <0$, so that the
fragmentation loses mass to dust and reaches $\mathbf{0}$ in finite time
a.s. As noticed in \cite{HaasMiermont}, the genealogy of the fragmentation
can then be described in terms of a continuum random tree.

The definition of CRT we are considering here is the one given by Aldous %
\cite{Aldous CRT3}, to which we refer for background and precise
definitions. Let $l_{1}:=\{\mathbf{x}=(x_{1,}x_{2},...),\sum_{k\geq 1}\left|
x_{k}\right| <\infty \}$ be endowed with the norm $\left\| x\right\|
_{1}:=\sum_{k\geq 1}\left| x_{k}\right| $, and let $\{\mathbf{e}_{k},k\geq
1\}$ be its usual basis. Roughly, a CRT is a pair $\left( \mathcal{T},\mu
\right) $ where $\mathcal{T}$ is a closed subset of $l_{1}$ that possesses
the ``tree'' property: for all $v,w$ $\in \mathcal{T}$, there exists a
unique (injective) path connecting $v$ to $w$, denoted by $\left[ \left[ v,w%
\right] \right] $. This tree is \textit{rooted}, that is one vertex is
distinguished as being the root $\emptyset _{\mathcal{T}}$. It is moreover
equipped with a $\sigma $-finite mass measure $\mu ,$ which is non-atomic
and puts mass only on the set of leaves, a leaf of $\mathcal{T}$ being a
vertex that does not belong to $\left[ \left[ \emptyset ,v\right[ \right[ $,
$\forall v\in \mathcal{T}$.

According to Theorem 1 of \cite{HaasMiermont}, since $\alpha <0,$ the
fragmentation $F^{(1)}$ can be constructed from some random compact CRT $%
\left( \mathcal{T}^{1},\mu ^{1}\right) $ rooted at $\mathbf{0}$\textbf{\ }as
follows: for each $t\geq 0$, $F^{(1)}(t)$ is the decreasing rearrangement of
the $\mu ^{1}$-masses of connected components of $\mathcal{T}^{1}$ obtained
by removing the vertices with a distance from the root smaller than $t$. We
shall say that $\left( \mathcal{T}^{1},\mu ^{1}\right) $ \textit{codes} the
fragmentation $F^{(1)}$. Note that the measure $\mu ^{1}$ is here a (random)
probability measure.

Now, let $\mathcal{T}^{m}$ denote the tree $\mathcal{T}^{1}$ rescaled by a
factor $m^{-\alpha }$ and let $\mu ^{m}$ be $m$ times the image measure of $%
\mu ^{1}$ by this scaling. Then, according to the self-similarity property, $%
(\mathcal{T}^{m},\mu ^{m})$ codes an ($\alpha ,\nu $)-fragmentation $F^{(m)}$%
.

\bigskip

In the remainder of this section we assume that%
\begin{equation}
m^{\alpha }\nu _{m}\rightarrow I,\text{ }I(l_{1}^{\downarrow })\neq 0,\text{
as }m\rightarrow \infty .  \label{32}
\end{equation}%
Given Theorem \ref{Thprincipal}, one can then expect that the sequence of
CRTs $(\mathcal{T}^{m},\mu ^{m})$ converges in distribution to some ``($%
\alpha ,\nu ,I$) fragmentation with immigration CRT'' ($\mathcal{T}_{FI},\mu
_{FI}$), which should be seen as an infinite baseline $\mathcal{B}:=\left\{ x%
\mathbf{e}_{1},x\geq 0\right\} $ on which compact CRTs are branched. A
version of this tree with a spine is constructed below.

We first specify the notion of convergence of trees we are using here. Two
trees are considered to be equivalent if there exists an isometry that maps
one onto the other and that preserves the root. Implicitly, we always
identify a tree with its equivalence class. A natural distance to consider
then is the so-called \textit{Gromov-Hausdorff} distance, which is a
distance measuring how far two metric spaces are from being isometric (see %
\cite{Gromov} for a precise definition and properties). Restricted to
compact trees of $l_{1}$, this distance is given by%
\begin{equation*}
D_{\mathcal{GH}}(\mathcal{T},\mathcal{T}^{\prime }):\mathcal{=}\inf (D_{%
\mathcal{H}}^{E}(\varphi (\mathcal{T}),\varphi ^{\prime }(\mathcal{T}%
^{\prime }))\vee d_{E}(\varphi (\emptyset _{\mathcal{T}}),\varphi ^{\prime
}(\emptyset _{\mathcal{T}^{\prime }})))
\end{equation*}%
where the infimum is taken over all isometric embeddings $\varphi ,\varphi
^{\prime }:l_{1}\rightarrow E$ into a same metric space $(E,d_{E})$ and $D_{%
\mathcal{H}}^{E}$ denotes the usual Hausdorff distance on the set of compact
subsets of $E$. However, the trees that appear as limit of $\mathcal{T}^{m}$
are not compact (but their restrictions to closed balls are). Hence we have
to trunk them, by introducing, for every tree $\mathcal{T}$ and every
integer $n,$ $\mathcal{T}\mid _{n}:=\{\mathbf{x\in }\mathcal{T}:\left\|
\mathbf{x}\right\| _{1}\leq n\}$. We then consider that a sequence $\mathcal{%
T}_{k}$ converges to $\mathcal{T}$ as $k\rightarrow \infty $ i.f.f. $D_{%
\mathcal{GH}}(\mathcal{T}_{k}\mid _{n},\mathcal{T}\mid _{n})\rightarrow 0$
for all $n\geq 0.$

\bigskip

Let us now construct a nice version of the CRT $(\mathcal{T}^{m},\mu ^{m})$
by using the geometric description of $F^{(m)}$ of Section 2.1. Instead of
branching fragmentations on a baseline, we here branch CRTs. To do so, write
$\mathbb{N}\backslash \left\{ 1\right\} =\biguplus_{i,j\geq 1}J_{i,j}$ where
Card$(J_{i,j})=\infty $ and let $f_{i,j}$ be a bijection between $\mathbb{N}$
and $J_{i,j}$. Remind then that $((r_{i},\mathbf{u}^{i}),i\geq 1)$ is a PPP
with intensity $I$ and that the random variables $r_{i}^{m}$, $\mathbf{u}%
^{i,m}$, $i\geq 1,$ introduced in formula (\ref{3}) Section 2.1, have been
constructed so that $\sum\nolimits_{i\geq 1}\delta _{(r_{i}^{m},\mathbf{u}%
^{i,m})}$ converges a.s. to $\sum\nolimits_{i\geq 1}\delta _{(r_{i},\mathbf{u%
}^{i})}$. Define then the maps%
\begin{eqnarray*}
m_{i,j}^{m} &:&\sum\nolimits_{k\geq 1}x_{k}\mathbf{e}_{k}\mapsto r_{i}^{m}%
\mathbf{e}_{1}+(u_{j}^{i,m})^{-\alpha }\sum\nolimits_{k\geq 1}x_{k}\mathbf{e}%
_{f_{i,j}(k)} \\
m_{i,j}^{I} &:&\sum\nolimits_{k\geq 1}x_{k}\mathbf{e}_{k}\mapsto r_{i}%
\mathbf{e}_{1}+(u_{j}^{i})^{-\alpha }\sum\nolimits_{k\geq 1}x_{k}\mathbf{e}%
_{f_{i,j}(k)}.
\end{eqnarray*}%
Introduce next a family $(\mathcal{T}_{i,j},\mu _{i,j}),i,j\geq 1,$ of
independent copies of $\left( \mathcal{T}^{1},\mu ^{1}\right) ,$ independent
of $(r_{i}^{m},\mathbf{u}^{i,m},i\geq 1)$ and $((r_{i},\mathbf{u}^{i}),i\geq
1).$ The tree $\mathcal{T}^{(u_{j}^{i,m})}:=m_{i,j}^{m}(\mathcal{T}_{i,j})$,
endowed with the measure $u_{j}^{i,m}\mu _{i,j}\circ (m_{i,j}^{m})^{-1}$,
codes a fragmentation $F^{(u_{j}^{i,m})}$ branched on $\mathcal{B}$ at
height $r_{i}^{m}$ and the required version of $(\mathcal{T}^{m},\mu ^{m})$,
still denoted by $(\mathcal{T}^{m},\mu ^{m}),$ is defined by
\begin{eqnarray}
\mathcal{T}^{m} &:&=\{x\mathbf{e}_{1},0\leq x\leq t_{\infty }^{m}\}\cup
_{i,j\geq 1}\mathcal{T}^{(u_{j}^{i,m})}  \label{21} \\
\mu ^{m} &:&=\sum\nolimits_{i,j\geq 1}u_{j}^{i,m}\mu _{i,j}\circ
(m_{i,j}^{m})^{-1},  \notag
\end{eqnarray}%
where $t_{\infty }^{m}$ is the first time at which $\Lambda ^{(m)}$ (defined
by (\ref{2})) reaches $0$.

Similarly, a nice version of the $(\alpha ,\nu ,I)$ fragmentation with
immigration CRT $(\mathcal{T}_{FI},\mu _{FI})$ is defined by%
\begin{eqnarray}
\mathcal{T}_{FI} &:&=\mathcal{B}\cup _{i,j\geq 1}\mathcal{T}^{(u_{j}^{i})}
\label{22} \\
\mu _{FI} &:&=\sum\nolimits_{i,j\geq 1}u_{j}^{i}\mu _{i,j}\circ
(m_{i,j}^{I})^{-1}  \notag
\end{eqnarray}%
where $\mathcal{T}^{(u_{j}^{i})}:=m_{i,j}^{I}(\mathcal{T}_{i,j})$. To obtain
a version of the $(\alpha ,\nu ,I)$ fragmentation with immigration from this
tree, just set $FI(t)$ for the decreasing sequence of $\mu _{FI}$-masses of
connected components of $\{\mathbf{x\in }\mathcal{T}_{FI}:\left\| \mathbf{x}%
\right\| _{1}\geq t,x_{1}\leq t\}$. At last, note that since $I$ is $%
(-\alpha )$-self-similar (by Lemma \ref{Lemma 6}), the CRT is also
self-similar, i.e.%
\begin{equation*}
(\mathcal{T}_{FI}^{a},a^{-1/\alpha }\mu _{FI}^{a})\overset{\text{law}}{=}%
\left( \mathcal{T}_{FI},\mu _{FI}\right) \text{ for all }a>0
\end{equation*}%
where $\mathcal{T}_{FI}^{a}$ is the tree $\mathcal{T}_{FI}$ rescaled by the
factor $a$ and $\mu _{FI}^{a}$ is the image measure of $\mu _{FI}$ by this
scaling.

We are now ready to state the counterpart, in term of trees, of Theorem \ref%
{Thprincipal}, assuming that (\ref{32}) holds. The topology on the set of
measures on $l_{1}$ is the topology of vague convergence.

\begin{theorem}
\label{PropTrees} As $m\rightarrow \infty $,%
\begin{equation*}
\left( \mathcal{T}^{m},\mu ^{m}\right) \overset{\mathrm{law}}{\rightarrow }%
\left( \mathcal{T}_{FI},\mu _{FI}\right) \text{.}
\end{equation*}
\end{theorem}

For the proof, we need the following lemma, where $h_{i,j}:=\sup \left\{
\left\| \mathbf{x}\right\| _{1},\mathbf{x}\in \mathcal{T}_{i,j}\right\} $ is
the height of the tree $\mathcal{T}_{i,j}$. It is known (see \cite{Haas1})
that those random variables have exponential moments.

\begin{lemma}
For all $n\in \mathbb{N}$,%
\begin{equation}
\sum\nolimits_{r_{i}^{m}\leq n,j\geq 1}u_{j}^{i,m}h_{i,j}^{-1/\alpha }%
\overset{\mathrm{P}}{\rightarrow }\sum\nolimits_{r_{i}\leq n,j\geq
1}u_{j}^{i}h_{i,j}^{-1/\alpha }\text{ \ \ as }m\rightarrow \infty \text{.}
\label{14}
\end{equation}%
As a consequence, one can extract from any increasing integer-valued
sequence $\boldsymbol{\kappa }$ a subsequence $\overline{\boldsymbol{\kappa }%
}$ such that for all $n,p\in \mathbb{N}$, as $m\rightarrow \infty $,%
\begin{equation}
\sum\nolimits_{r_{i}^{\overline{\boldsymbol{\kappa }}_{m}}\leq n,j\geq 1}%
\mathbf{1}_{\{(u_{j}^{i,\overline{\boldsymbol{\kappa }}_{m}})^{-\alpha
}h_{i,j}p>1\}}\overset{\mathrm{a.s.}}{\rightarrow }\sum\nolimits_{r_{i}\leq
n,j\geq 1}\mathbf{1}_{\{(u_{j}^{i})^{-\alpha }h_{i,j}p>1\}}<\infty .
\label{15}
\end{equation}
\end{lemma}

\noindent \textbf{Proof.} \textbf{(i).} Fix $n\in \mathbb{N}$ and recall
that a.s. $\sum\nolimits_{i\geq 1}\delta _{(r_{i}^{m},\mathbf{u}^{i,m})}$
converges to $\sum\nolimits_{i\geq 1}\delta _{(r_{i},\mathbf{u}^{i})}$, and $%
r_{i}\notin \mathbb{N}$, $i\geq 1.$ Consequently, $u_{j}^{i,m}\mathbf{1}%
_{\{r_{i}^{m}\leq n\}}\rightarrow u_{j}^{i}\mathbf{1}_{\{r_{i}\leq n\}}$ for
all $i,j\geq 1$ a.s., and a.s. for all $\eta >0$, there exists a $k\in
\mathbb{N}$ such that for $m$ large enough,
\begin{equation}
\sum\nolimits_{i+j\geq k}(u_{j}^{i,m}\mathbf{1}_{\{r_{i}^{m}\leq
n\}}+u_{j}^{i}\mathbf{1}_{\{r_{i}\leq n\}})\leq \eta .  \label{1}
\end{equation}%
We want to prove that $X_{m}:=\sum\nolimits_{r_{i}^{m}\leq n,j\geq
1}u_{j}^{i,m}h_{i,j}^{-1/\alpha }$ converges to $X:=\sum\nolimits_{r_{i}\leq
n,j\geq 1}u_{j}^{i}h_{i,j}^{-1/\alpha }$ in probability. Remark that $%
X<\infty $ a.s. since $E[X\mid (r_{i},\mathbf{u}^{i}),i\geq
1]=E[h_{1,1}^{-1/\alpha }]\sum\nolimits_{r_{i}\leq n,j\geq 1}u_{j}^{i}$ is
finite a.s. Similarly, $X_{m}<\infty $ a.s. Then, since%
\begin{equation*}
P\left( \left| X_{m}-X\right| >\varepsilon \right) =E\left[ E\left[ \mathbf{1%
}_{\{\left| X_{m}-X\right| >\varepsilon \}}\mid (r_{i}^{m},\mathbf{u}%
^{i,m}),(r_{i},\mathbf{u}^{i}),i,m\geq 1\right] \right] ,
\end{equation*}%
it is sufficient, by dominated convergence, to prove that the conditional
expectation converges a.s. to $0$, $\forall \varepsilon >0.$ For large $m$%
's, one has%
\begin{eqnarray*}
E[\mathbf{1}_{\{\left| X_{m}-X\right| >\varepsilon \}} &\mid &(r_{i}^{m},%
\mathbf{u}^{i,m}),(r_{i},\mathbf{u}^{i}),i,m\geq 1] \\
&\leq &\varepsilon ^{-1}E\left[ \left| X_{m}-X\right| \mid (r_{i}^{m},%
\mathbf{u}^{i,m}),(r_{i},\mathbf{u}^{i}),i,m\geq 1\right]  \\
&\leq &\varepsilon ^{-1}E[h_{1,1}^{-1/\alpha }]\sum\nolimits_{i,j\geq
1}\left| u_{j}^{i,m}\mathbf{1}_{\{r_{i}^{m}\leq n\}}-u_{j}^{i}\mathbf{1}%
_{\{r_{i}\leq n\}}\right|  \\
&\leq &\varepsilon ^{-1}E[h_{1,1}^{-1/\alpha }]\left(
\sum\nolimits_{i+j<k}\left| u_{j}^{i,m}\mathbf{1}_{\{r_{i}^{m}\leq
n\}}-u_{j}^{i}\mathbf{1}_{\{r_{i}\leq n\}}\right| +\eta \right) ,
\end{eqnarray*}%
the last inequality coming from (\ref{1}). So for all $\eta >0$, we have a
upper bound smaller than $2\eta \varepsilon ^{-1}E[h_{1,1}^{-1/\alpha }]$
for all $m$ sufficiently large, a.s. Hence the conclusion.

\textbf{(ii).} The measure $I$ is self-similar (by Lemma \ref{Lemma 6}) and
consequently atomless on $l_{1}^{\downarrow }\backslash \{\mathbf{0}\}$. As $%
(r_{i},\mathbf{u}^{i})_{i\geq 1}$ is a PPP with intensity $I$, independent
of the $h_{i,j}$'s, this implies that a.s. $(u_{j}^{i})^{-\alpha
}h_{i,j}p\neq 1$, $\forall i,j,p\geq 1,$ which in turn leads to the
convergence of $\mathbf{1}_{\{(u_{j}^{i,m})^{-\alpha }h_{i,j}p>1\}}\mathbf{1}%
_{\{r_{i}^{m}\leq n\}}$\ to $\mathbf{1}_{\{(u_{j}^{i})^{-\alpha
}h_{i,j}p>1\}}\mathbf{1}_{\{r_{i}\leq n\}}$ a.s. $\forall i,j,p,n\geq 1$.
Then for all $k\geq 1,$%
\begin{eqnarray}
&&\left| \sum\nolimits_{r_{i}^{m}\leq n,j\geq 1}\mathbf{1}%
_{\{(u_{j}^{i,m})^{-\alpha }h_{i,j}p>1\}}-\sum\nolimits_{r_{i}\leq n,j\geq 1}%
\mathbf{1}_{\{(u_{j}^{i})^{-\alpha }h_{i,j}p>1\}}\right|  \label{33} \\
&\leq &\sum\nolimits_{i+j<k}\left| \mathbf{1}_{\{(u_{j}^{i,m})^{-\alpha
}h_{i,j}p>1\}}\mathbf{1}_{\{r_{i}^{m}\leq n\}}-\mathbf{1}_{\{(u_{j}^{i})^{-%
\alpha }h_{i,j}p>1\}}\mathbf{1}_{\{r_{i}\leq n\}}\right|  \notag \\
&&+p^{-1/\alpha }\sum\nolimits_{i+j\geq k}(u_{j}^{i,m}h_{i,j}^{-1/\alpha }%
\mathbf{1}_{\{r_{i}^{m}\leq n\}}+u_{j}^{i}h_{i,j}^{-1/\alpha }\mathbf{1}%
_{\{r_{i}\leq n\}}).  \notag
\end{eqnarray}%
So if we prove that each sequence $\boldsymbol{\kappa }$ possesses a
subsequence $\overline{\boldsymbol{\kappa }}$ \textit{independent of} $n\in
\mathbb{N}$ such that, a.s. for all $\varepsilon >0$ there exists a $k$ such
that%
\begin{equation}
\sum\nolimits_{i+j\geq k}(u_{j}^{i,\overline{\boldsymbol{\kappa }}%
_{m}}h_{i,j}^{-1/\alpha }\mathbf{1}_{\{r_{i}^{\overline{\boldsymbol{\kappa }}%
_{m}}\leq n\}}+u_{j}^{i}h_{i,j}^{-1/\alpha }\mathbf{1}_{\{r_{i}\leq
n\}})\leq \varepsilon \text{ for all }m\text{ large enough,}  \label{18}
\end{equation}%
then we will have the statement (using also that the first term in the right
hand side of the inequality (\ref{33}) is composed by a finite number of
terms that all converge to $0$ as $m\rightarrow \infty $). Clearly, to get (%
\ref{18}), it is sufficient to show that there is a subsequence $\overline{%
\boldsymbol{\kappa
}}$ such that $\forall n$,
\begin{equation*}
\sum\nolimits_{r_{i}^{\overline{\boldsymbol{\kappa }}_{m}}\leq n,j\geq
1}u_{j}^{i,\overline{\boldsymbol{\kappa }}_{m}}h_{i,j}^{-1/\alpha
}\rightarrow \sum\nolimits_{r_{i}\leq n,j\geq 1}u_{j}^{i}h_{i,j}^{-1/\alpha }%
\text{ a.s.}
\end{equation*}%
To construct this subsequence, we use the convergence in probability (\ref%
{14}). It implies that for all $n$, there is a subsequence $\overline{%
\boldsymbol{\kappa }}^{(n)}$ such that the above a.s. convergence holds. We
want a sequence $\overline{\boldsymbol{\kappa }}$ independent of $n$ and to
do so, use a diagonal extraction argument: extract $\overline{%
\boldsymbol{\kappa }}^{(1)}$ from $\boldsymbol{\kappa }$ and then
recursively $\overline{\boldsymbol{\kappa }}^{(n+1)}$ from $\overline{%
\boldsymbol{\kappa }}^{(n)}$. Then set $\overline{\boldsymbol{\kappa }}_{m}:=%
\overline{\boldsymbol{\kappa }}^{(m)}(m).$\ \ \rule{0.5em}{0.5em}

\textbf{\newline
Proof of Theorem \ref{PropTrees}.} In all the proof ($\mathcal{T}^{m},\mu
^{m}$) refers to the version (\ref{21}) of the fragmentation CRT with total
weight $m$ and ($\mathcal{T}_{FI},\mu _{FI}$) to the version (\ref{22}) of
the fragmentation with immigration CRT. We will prove that ($\mathcal{T}%
^{m},\mu ^{m}$) converges in probability\ to ($\mathcal{T}_{FI},\mu _{FI}$),
or equivalently that for any increasing integer-valued sequence $%
\boldsymbol{\kappa }$, one can extract a subsequence $\overline{%
\boldsymbol{\kappa }}$ such that ($\mathcal{T}^{\overline{%
\boldsymbol{\kappa
}}_{m}},\mu ^{\overline{\boldsymbol{\kappa }}_{m}}$) converges a.s. to ($%
\mathcal{T}_{FI},\mu _{FI}$). So, fix such a sequence $\boldsymbol{\kappa
}$ and consider its subsequence $\overline{\boldsymbol{\kappa }}$ introduced
in the lemma above, so that the a.s. convergences (\ref{15}) hold. In the
rest of the proof, all the assertions hold a.s., so we drop the ``a.s.''
from the notations.

\textbf{(i).} A first remark is that for all $i,j\geq 1$,
\begin{equation}
D_{\mathcal{H}}^{l_{1}}(\mathcal{T}^{(u_{j}^{i,m})},\mathcal{T}%
^{(u_{j}^{i})})\leq \left| r_{i}^{m}-r_{i}\right| +\left|
(u_{j}^{i,m})^{-\alpha }-(u_{j}^{i})^{-\alpha }\right| h_{i,j}\rightarrow 0%
\text{ \ as }m\rightarrow \infty .  \label{17}
\end{equation}%
Fix then $n,p\in \mathbb{N}$. As a consequence of (\ref{15}), the number of
trees among $\{\mathcal{T}^{(u_{j}^{i,\overline{\boldsymbol{\kappa }}_{m}})}$%
, $i,j\geq 1$, $r_{i}^{\overline{\boldsymbol{\kappa }}_{m}}\leq n\}$, which
are not entirely contained in $\{\mathbf{x}:\left\| \mathbf{x-}x_{1}\mathbf{e%
}_{1}\right\| _{1}\leq p^{-1}\}$ is constant (finite) for $m$ large enough.
Let $\mathcal{K}$ be the finite set of $(i,j)$ s.t. $\mathcal{T}%
^{(u_{j}^{i})}$ is not entirely contained in $\{\mathbf{x}:\left\| \mathbf{x-%
}x_{1}\mathbf{e}_{1}\right\| _{1}\leq p^{-1}\}$. Then for large $m$'s,
\begin{equation*}
D_{\mathcal{H}}^{l_{1}}(\mathcal{T}^{\overline{\boldsymbol{\kappa }}%
_{m}}\mid _{n},\mathcal{T}_{FI}\mid _{n})\leq p^{-1}+\max_{i,j\in \mathcal{K}%
}D_{\mathcal{H}}^{l_{1}}(\mathcal{T}^{(u_{j}^{i,\overline{\boldsymbol{\kappa
}}_{m}})},\mathcal{T}^{(u_{j}^{i})}).
\end{equation*}%
Considering (\ref{17}) and taking $m$ larger if necessary, one sees that
this upper bound is in turn bounded by $2p^{-1}$. This holds for all $p\in
\mathbb{N}$. Hence $D_{\mathcal{GH}}(\mathcal{T}^{\overline{%
\boldsymbol{\kappa }}_{m}}\mid _{n},\mathcal{T}_{FI}\mid _{n})\rightarrow 0$%
, $\forall n\in \mathbb{N}$.

\textbf{(ii).} Next, for all $\mathbb{R}^{+}$-valued continuous function $f$
with compact support on $l_{1}$, \linebreak $\langle u_{j}^{i,m}\mu
_{i,j}\circ (m_{i,j}^{m})^{-1},f\rangle $ converges to $\langle u_{j}^{i}\mu
_{i,j}\circ (m_{i,j}^{I})^{-1},f\rangle $, since $m_{i,j}^{m}(\mathbf{x}%
)\rightarrow m_{i,j}^{I}(\mathbf{x})$ for all $\mathbf{x}\in l_{1}$. To
deduce from this that the sum over $i,j\geq 1$ of these measures converges,
fix some $\eta >0$ and let $C_{f}:=\sup_{l_{1}}\left| f(\mathbf{x})\right| $%
.\ Again we use the argument that there exists some $k\in \mathbb{N}$ such
that $\sum\nolimits_{i+j\geq k}u_{j}^{i}\mathbf{1}_{\{r_{i}\leq
C_{f}\}}<\eta $\ and $\sum\nolimits_{i+j\geq k}u_{j}^{i,m}\mathbf{1}%
_{\{r_{i}^{m}\leq C_{f}\}}<\eta $\ for all $m$ large enough, which leads to
\begin{equation*}
\left| \left\langle \mu ^{m},f\right\rangle -\left\langle \mu
_{FI},f\right\rangle \right| \leq 2C_{f}\eta +\sum\nolimits_{i+j<k}\left|
\langle u_{j}^{i,m}\mu _{i,j}\circ (m_{i,j}^{m})^{-1},f\rangle -\langle
u_{j}^{i}\mu _{i,j}\circ (m_{i,j}^{I})^{-1},f\rangle \right|
\end{equation*}%
which is bounded by $(2C_{f}+1)\eta $ for large $m$'s$.$ \ \rule%
{0.5em}{0.5em}

\section{Stable fragmentations}

In this section, we apply our results to two specific families of
fragmentations constructed from the so-called stable tree $(\mathcal{T}%
^{\beta },\mu ^{\beta })$ with index $\beta $, $1<\beta <2$. This object is
a CRT introduced by Duquesne and Le Gall \cite{DuquesneLegall}, who we refer
to for a rigorous construction. Roughly, $\mathcal{T}^{\beta }$ arises as
the limit in distribution of rescaled critical Galton-Watson trees $T_{n}$,
conditioned to have $n$ vertices and edge-lengths $n^{\beta ^{-1}-1}$, and
an offspring distribution $(\eta _{k},k\geq 0)$ such that $\eta _{k}\sim
Ck^{-1-\beta }$ as $k\rightarrow \infty $. It is endowed with a probability
measure $\mu ^{\beta }$ which is the limit as $n\rightarrow \infty $ of the
empiric measure on the vertices of $T_{n}$.

\subsection{Stable fragmentations with a negative index of self-similarity}

Let $F^{\beta -}(t)$ denotes the decreasing sequence of $\mu ^{\beta }$%
-masses of connected components obtained by removing in $\mathcal{T}^{\beta
} $ all vertices at distance less than $t$ from the root, $t\geq 0$.
Miermont \cite{mierfmoins} shows that $F^{\beta -}$ is a self-similar
fragmentation with index $1/\beta -1$, and with a dislocation measure $\nu
^{\beta }$ given by%
\begin{equation*}
\int_{l_{1,\leq 1}^{\downarrow }}f(\mathbf{s})\nu ^{\beta }(\mathrm{d}%
\mathbf{s})=C_{\beta }E\left[ T_{1}^{\beta }f((T_{1}^{\beta })^{-1}(\Delta
_{1}^{\beta },\Delta _{2}^{\beta },...))\right] \text{, }f\in \mathcal{F},
\end{equation*}%
where $C_{\beta }=\beta ^{2}\Gamma (2-\beta ^{-1})/\Gamma \left( 2-\beta
\right) .$ The process $T^{\beta }$ is a stable subordinator with Laplace
exponent $q^{1/\beta }$, i.e.
\begin{equation}
E[\exp (-qT_{r}^{\beta })]=\exp (-rq^{1/\beta })\text{, }q,r\geq 0\text{,}
\label{23}
\end{equation}%
and $(\Delta _{1}^{\beta },\Delta _{2}^{\beta },...)$ denotes the sequence
of jumps in the decreasing order of $T^{\beta }$ before time $1.$ In order
to apply Theorem \ref{Thprincipal} to these fragmentations, we state the
following lemma.

\begin{lemma}
\label{Lemmastable}As $m\rightarrow \infty $, $m^{1/\beta -1}\nu _{m}^{\beta
}\rightarrow I^{\beta }$, where $I^{\beta }$ is defined by%
\begin{equation*}
\int_{l_{1}^{\downarrow }}f(\mathbf{s})I^{\beta }(\mathrm{d}\mathbf{s}%
)=\beta (\beta -1)(\Gamma (2-\beta ))^{-1}\int_{0}^{\infty }E\left[
f(x^{\beta }(\Delta _{1}^{\beta },\Delta _{2}^{\beta },...))\right]
x^{-\beta }\mathrm{d}x,\text{ }f\in \mathcal{F}.
\end{equation*}
\end{lemma}

Using (\ref{23}), one sees that $I^{\beta }$ integrates $(1-\exp
(-\sum\nolimits_{i\geq 1}s_{i}))$ and therefore that it is an immigration
measure (it is also a consequence of the above convergence).

\textbf{\newline
Proof.} In all the proof, $T_{1}^{\beta },\Delta _{1}^{\beta },\Delta
_{2}^{\beta },...$ are rather denoted by $T_{1},\Delta _{1},\Delta _{2},...$
A classical idea is to use a size-biased permutation $(\Delta _{1}^{\ast
},\Delta _{2}^{\ast },...)$ of $(\Delta _{1},\Delta _{2},...)$ to obtain
some results on the latter. To do so, we first recall that $T_{1}$ has a
density (see e.g. formula (40) in \cite{Pitman}), that we denote by $q$. One
then obtains, using Palm measures theory (see e.g. \cite{Neveu}), the
following equality :%
\begin{eqnarray}
&&E\left[ f(T_{1},\Delta _{1}^{\ast },\Delta _{2}^{\ast },\Delta _{3}^{\ast
},...,\Delta _{k+1}^{\ast })\right]  \label{7} \\
&=&c_{\beta }^{k+1}\int_{0}^{\infty
}\int_{0}^{s_{0}}\int_{0}^{s_{0-}s_{1}}...\int_{0}^{s_{0}-s_{1}-...-s_{k}}%
\frac{f(s_{0},s_{1},...,s_{k+1})q(s_{0}-s_{1}-...-s_{k+1})\mathrm{d}%
s_{k+1}...\mathrm{d}s_{1}\mathrm{d}s_{0}}{s_{k+1}^{1/\beta }s_{k}^{1/\beta
}...s_{1}^{1/\beta }(s_{0}-s_{1}-...-s_{k})...(s_{0}-s_{1})s_{0}}  \notag
\end{eqnarray}%
for all non-negative measurable function $f$ on $(\mathbb{R}^{+})^{k+2}$,
where $c_{\beta }=(\beta \Gamma (1-1/\beta ))^{-1}$.

Let then $g$ be a non-negative measurable function on $(\mathbb{R}^{+})^{k}$%
. One has%
\begin{eqnarray}
&&m^{1/\beta -1}E\left[ T_{1}g\left( \frac{m\Delta _{2}^{\ast }}{T_{1}},%
\frac{m\Delta _{3}^{\ast }}{T_{1}},...,\frac{m\Delta _{k+1}^{\ast }}{T_{1}}%
\right) \right]  \label{25} \\
&=&c_{\beta }^{k+1}\int_{0}^{\infty
}\int_{0}^{u}\int_{0}^{u-s_{2}}...\int_{0}^{u-...-s_{k}}\left( m^{1/\beta
-1}\int_{u}^{\infty }\frac{g(ms_{2}/s_{0},ms_{3}/s_{0},...,ms_{k+1}/s_{0})}{%
(s_{0}-u)^{1/\beta }}\mathrm{d}s_{0}\right)  \notag \\
&&\times \frac{q(u-s_{2}...-s_{k+1})\mathrm{d}s_{k+1}...\mathrm{d}s_{2}%
\mathrm{d}u}{s_{k+1}^{1/\beta }s_{k}^{1/\beta }...s_{2}^{1/\beta
}(u-s_{2}...-s_{k})(u-s_{2}...-s_{k-1})...(u-s_{2})u}  \notag \\
&=&\beta c_{\beta }E\left[ \int_{0}^{\infty }\frac{g(\Delta _{1}^{\ast
}/(v^{-\beta }+m^{-1}T_{1}),...,\Delta _{k}^{\ast }/(v^{-\beta
}+m^{-1}T_{1}))}{v^{\beta }}\mathrm{d}v\right]  \notag
\end{eqnarray}%
where for the first equality we use formula (\ref{7}), the change of
variables $s_{1}\mapsto s_{0}-u$ and Fubini's Theorem, and for the second
equality the change of variables $s_{0}\mapsto u+mv^{-\beta }$ and again
formula (\ref{7}). This holds in particular for $g(x_{1},...,x_{k})=f\circ
d^{\downarrow }(x_{1},...,x_{k},0,...)$ when $f\in \mathcal{F}$ and $%
d^{\downarrow }$ is the function that associates to $(x_{1},x_{2},...)\in (%
\mathbb{R}^{+})^{\infty },$ $\sum_{i\geq 1}x_{i}<\infty ,$ its decreasing
rearrangement in $l_{1}^{\downarrow }$ (this function is measurable). Our
aim now is to let $k\rightarrow \infty $ in equality (\ref{25}) for such
functions $g$.\ To do so, first note that $d^{\downarrow
}(x_{1},...,x_{k},0,...)\rightarrow d^{\downarrow }(x_{1},x_{2},...)$ in $%
l_{1}^{\downarrow }$ as $k\rightarrow \infty $, for all $(x_{1},x_{2},...)%
\in (\mathbb{R}^{+})^{\infty },$ $\sum_{i\geq 1}x_{i}<\infty $. We then
claim that dominated convergence applies in both sides of the equality.
Indeed, for the left hand side, since $f(\mathbf{s})\leq
(\sum\nolimits_{i\geq 1}s_{i})\wedge 1$, one has for all $k$,
\begin{equation*}
T_{1}g\left( \frac{m\Delta _{2}^{\ast }}{T_{1}},\frac{m\Delta _{3}^{\ast }}{%
T_{1}},...,\frac{m\Delta _{k+1}^{\ast }}{T_{1}}\right) \leq T_{1}\left(
1\wedge m\left( \frac{T_{1}-\Delta _{1}^{\ast }}{T_{1}}\right) \right) .
\end{equation*}%
It is therefore sufficient to prove that $E[T_{1}\wedge m(T_{1}-\Delta
_{1}^{\ast })]<\infty $, which, clearly, holds if $E[T_{1}-\Delta _{1}^{\ast
}]<\infty $. It is not hard to see, using the joint distribution (\ref{7}),
that the last expectation is bounded from above (up to a finite constant) by
$E[(T_{1})^{1-1/\beta }]$, which, according to formula (43) in \cite{Pitman}%
, is finite. Hence dominated convergence applies in the left hand side of (%
\ref{25}). Now, for the right hand side, one uses that%
\begin{equation*}
g(\Delta _{1}^{\ast }/(v^{-\beta }+m^{-1}T_{1}),...,\Delta _{k}^{\ast
}/(v^{-\beta }+m^{-1}T_{1}))\leq (T_{1}v^{\beta }\wedge 1)
\end{equation*}%
which is integrable with respect to d$\mathbb{P\otimes }v^{-\beta }$d$v$,
because $(1-\exp (-T_{1}v^{\beta }))$ is. At last, letting $k\rightarrow
\infty $, one obtains%
\begin{equation*}
\begin{array}{ll}
m^{1/\beta -1}E\left[ T_{1}f\circ d^{\downarrow }\left( \frac{m\Delta
_{2}^{\ast }}{T_{1}},\frac{m\Delta _{3}^{\ast }}{T_{1}},...\right) \right] &
=\beta c_{\beta }E\left[ \int_{0}^{\infty }\frac{f\circ d^{\downarrow
}(\Delta _{1}^{\ast }/(v^{-\beta }+m^{-1}T_{1}),\Delta _{2}^{\ast
}/(v^{-\beta }+m^{-1}T_{1}),...)}{v^{\beta }}\mathrm{d}v\right] \\
& =\beta c_{\beta }E\left[ \int_{0}^{\infty }\frac{f(\Delta _{1}/(v^{-\beta
}+m^{-1}T_{1}),\Delta _{2}/(v^{-\beta }+m^{-1}T_{1}),...)}{v^{\beta }}%
\mathrm{d}v\right] .%
\end{array}%
\end{equation*}%
This latter term converges as $m\rightarrow \infty $ to $\beta c_{\beta
}\int_{0}^{\infty }E\left[ f(v^{\beta }(\Delta _{1},\Delta _{2},...))\right]
v^{-\beta }\mathrm{d}v$, again by dominated convergence. Hence we would have
the required convergence $\langle m^{1/\beta -1}\nu _{m}^{\beta },f\rangle
\rightarrow \langle I^{\beta },f\rangle $ for all continuous non-negative
functions $f\in \mathcal{F}$ if we could replace in the left hand side of
the above formula the sequence $d^{\downarrow }(m\Delta _{2}^{\ast
}/T_{1},m\Delta _{3}^{\ast }/T_{1},...)$ by $(m\Delta _{2}/T_{1},m\Delta
_{3}/T_{1},...)$. Of course, this is not possible. However, conditional on $%
\Delta _{1}^{\ast }>T_{1}/2$, one has $\Delta _{1}^{\ast }\overset{\mathrm{%
a.s.}}{=}\Delta _{1}$ (equivalently $d^{\downarrow }\left( \Delta _{2}^{\ast
},\Delta _{3}^{\ast },...\right) \overset{\mathrm{a.s.}}{=}\left( \Delta
_{2},\Delta _{3},...\right) $), since the size-biased pick $\Delta
_{1}^{\ast }/T_{1}$ is then necessarily equal to the largest mass $\Delta
_{1}/T_{1}$. Therefore, one can write%
\begin{equation*}
\begin{array}{l}
m^{1/\beta -1}E\left[ T_{1}f\left( \frac{m\Delta _{2}}{T_{1}},\frac{m\Delta
_{3}}{T_{1}},...\right) \right] =m^{1/\beta -1}E\left[ T_{1}f\circ
d^{\downarrow }\left( \frac{m\Delta _{2}^{\ast }}{T_{1}},\frac{m\Delta
_{3}^{\ast }}{T_{1}},...\right) \right] \\
\text{ \ \ }+m^{1/\beta -1}E\left[ \left( T_{1}f\left( \frac{m\Delta _{2}}{%
T_{1}},\frac{m\Delta _{3}}{T_{1}},...\right) -T_{1}f\circ d^{\downarrow
}\left( \frac{m\Delta _{2}^{\ast }}{T_{1}},\frac{m\Delta _{3}^{\ast }}{T_{1}}%
,...\right) \right) \mathbf{1}_{\{\Delta _{1}^{\ast }\leq T_{1}/2\}}\right] .%
\end{array}%
\end{equation*}%
and this converges to the required limit as $m\rightarrow \infty $, because
the absolute value of the second term in the right hand side of the equality
is bounded from above by $m^{1/\beta -1}E\left[ 2T_{1}\mathbf{1}_{\{\Delta
_{1}^{\ast }\leq T_{1}/2\}}\right] $ which in turn is bounded by $m^{1/\beta
-1}E\left[ 4(T_{1}-\Delta _{1}^{\ast })\right] ,$ which converges to $0$ as $%
m\rightarrow \infty .$ \ \rule{0.5em}{0.5em}

\bigskip

From this and Theorem \ref{Thprincipal}, one deduces that

\begin{equation}
\left( m-(F_{1}^{\beta -})^{(m)},((F_{2}^{\beta -})^{(m)},(F_{3}^{\beta
-})^{(m)},...)\right) \overset{\mathrm{law}}{\rightarrow }\left( \sigma
_{I^{\beta }},FI^{\beta }\right)  \label{24}
\end{equation}%
where $FI^{\beta }$ is a fragmentation with immigration process $(1/\beta
-1,\nu ^{\beta },I^{\beta })$ and $\sigma _{I^{\beta }}$ is the stable
subordinator with index $1-1/\beta $ representing the total mass of
immigrants. In terms of trees (Theorem \ref{PropTrees}), one has%
\begin{equation*}
\left( (\mathcal{T}^{\beta ,m},m\mu ^{\beta ,m}\right) \overset{\mathrm{law}}%
{\rightarrow }\left( \mathcal{T}_{FI^{\beta }},\mu _{FI^{\beta }}\right)
\end{equation*}%
where $\mathcal{T}^{\beta ,m}$ is the stable tree rescaled by a factor $%
m^{1-1/\beta }$ and $\mu ^{\beta ,m}$ is the image of $\mu ^{\beta }$ by
this scaling; $\left( \mathcal{T}_{FI^{\beta }},\mu _{FI^{\beta }}\right) $
is a fragmentation with immigration CRT with parameters $\left( 1/\beta
-1,\nu ^{\beta },FI^{\beta }\right) $.

In chapter 4.4.2 of \cite{HaasThese}, it is shown that (some version of)
this fragmentation with immigration $FI^{\beta }$ can be constructed from
the height process $H^{\beta }$ coding a continuous state branching process
with immigration, with branching mechanism $u^{\beta }$ and immigration
mechanism $\beta u^{\beta -1}$ as follows: for all $t\geq 0,$ $FI^{\beta
}(t) $ is the decreasing rearrangement of the lengths of finite excursions
of $H^{\beta }$ above $t$. In a recent work, Duquesne \cite%
{DuquesneImmigration} shows that the rescaled height function of some
ordered version of the stable tree with index $\beta $ converges to $%
H^{\beta }$, which corroborates our results.

Last, thanks to the self-similarity, the convergence (\ref{24}) also
specifies the behavior of $(F^{\beta -})^{(1)}(\varepsilon \cdot )$ as $%
\varepsilon \rightarrow 0$. In particular, the mass of dust $1-(M^{\beta
})^{(1)}$ behaves as follows.

\begin{corollary}
As $\varepsilon \rightarrow 0$, $\varepsilon ^{-\beta /(\beta
-1)}(1-(M^{\beta })^{(1)}(\varepsilon \cdot ))\overset{\mathrm{law}}{%
\rightarrow }\int_{0}^{t}L^{\beta }(u)$\textrm{d}$u$, where $L^{\beta }$ is
a continuous state branching process with immigration starting from $0,$
with branching mechanism $u^{\beta }$ and immigration mechanism $\beta
u^{\beta -1}$.
\end{corollary}

Indeed, according to (\ref{24}), $\varepsilon ^{-\beta /(\beta
-1)}(1-(M^{\beta })^{(1)}(\varepsilon \cdot ))$ converges in law to some
non-trivial limit that corresponds to the total mass of microscopic
particles produced until time $t$ by the fragmentation with immigration $%
FI^{\beta }$ and the construction of $FI^{\beta }$ from $H^{\beta }$ implies
that this limit is equal to $\int_{0}^{t}L^{\beta }(u)$\textrm{d}$u$ where $%
L^{\beta }(u)$ is the local time at $u$ of $H^{\beta }.$ Lambert \cite%
{Lambert} proves that $L^{\beta }$ is actually a continuous state branching
process with immigration starting from $0,$ with the characteristics
mentioned in the above corollary. In a previous work, Miermont \cite%
{mierfmoins} obtained this convergence result on the mass of dust in the one
dimensional case.

\bigskip

\noindent \textbf{Remark. }Using the same tools, one sees that the above
corollary is also valid when replacing $\beta $ by $2$ and the fragmentation
$F^{\beta }$ by a self-similar fragmentation with index $-1/2$ and
dislocation measure $\sqrt{2}\nu _{B_{r}}$.

\subsection{Stable fragmentations with a positive index of self-similarity}

We here consider the fragmentations $F^{\beta +}$ with parameters $(1/\beta
,\nu ^{\beta })$, $1<\beta <2$. Such fragmentations can also be constructed
from the stable trees $\mathcal{T}^{\beta }$, by cutting them at nodes (see %
\cite{mierfplus}). According to Lemma \ref{Lemmastable} and Corollary \ref%
{Corosmallbeh}, one knows that their small times behaviors are characterized
in terms of the pure immigration processes $(I^{\beta }(t),t\geq 0)$ with
intensity $I^{\beta }$. Let then $\varrho $ be a $(\beta -1)$-stable
subordinator with Laplace exponent $\beta q^{\beta -1}$, $q\geq 0$,
independent of $T^{\beta }$ and call $(\Delta _{1}^{\beta }(\varrho
(t)),\Delta _{2}^{\beta }(\varrho (t)),...)$ the decreasing sequence of
jumps of $T^{\beta }$ before time $\rho (t),$ $t\geq 0.$ A moment of thought
shows that $((\Delta _{1}^{\beta }(\varrho (t)),\Delta _{2}^{\beta }(\varrho
(t)),...),t\geq 0)$ is distributed as $(I^{\beta }(t),t\geq 0)$. Therefore,

\begin{corollary}
As $\varepsilon \rightarrow 0$,%
\begin{equation*}
\varepsilon ^{-\beta /(\beta -1)}(1-F_{1}^{\beta +}(\varepsilon \cdot
),(F_{2}^{\beta +}(\varepsilon \cdot ),F_{3}^{\beta +}(\varepsilon \cdot
),...))\overset{\mathrm{law}}{\rightarrow }(T_{\varrho (\cdot )}^{\beta
},(\Delta _{1}^{\beta }(\varrho (\cdot )),\Delta _{2}^{\beta }(\varrho
(\cdot )),...)).
\end{equation*}
\end{corollary}

Miermont \cite{mierfplus} obtained this result in the one dimensional case.

\section{Appendix}

\subsection{Proof of Proposition \ref{Propmassecontinue}}

Our aim is to prove that under the general assumptions we have made on $\tau
,\nu $ ($\tau $ monotone near $0$, $\nu (\sum\nolimits_{i\geq 1}s_{i}<1)=0$)
the mass $M^{(m)}(t)=\sum\nolimits_{i\geq 1}F_{i}^{(m)}(t)$ is continuous in
$t$. The proof is the same for all $m,$ so we suppose that $m=1$ and we use
the notations $M,$ $F$ instead of $M^{(1)}$, $F^{(1)}$. We also suppose that
$\nu (l_{1,\leq 1}^{\downarrow })\neq 0$.

As often in the study of loss of mass, the problem can be tackled by
considering the evolution of some fragments independently tagged at random.
So, consider the interval representation $I^{\tau }$ from which $F$ has been
constructed in Section 1.1.1 and let $U,U^{\prime }$ be two independent r.v.
uniformly distributed on $\left( 0,1\right) $, independent of $I^{\tau }$.
Let then $D_{\tau }$ (resp. $D_{\tau }^{\prime }$) be the first time,
possibly infinite, at which $U$ (resp. $U^{\prime }$) falls into the dust
and note that with probability one, $P(D_{\tau }=D_{\tau }^{\prime }=t\mid
I^{\tau })=(M(t-)-M(t))^{2}$ for all $t\geq 0.$ Consequently, the mass $M$
is a.s. continuous as soon as $P(D_{\tau }=D_{\tau }^{\prime }<\infty )=0$.

The goal now is to prove that this probability is equal to $0$. To do so,
note first, using the time changes (\ref{4}), that
\begin{equation*}
D_{\tau }=\int_{0}^{\infty }\mathrm{d}r/\tau (I_{U}^{\hom
}(r))=\int_{0}^{\infty }\mathrm{d}r/\tau (\exp (-\sigma (r))),
\end{equation*}%
where, by definition, $\sigma =-\ln (I_{U}^{\hom })$. A well-known result of %
\cite{berthfrag01} says that $\sigma $ is a subordinator with \textit{zero
drift} and L\'{e}vy measure $L($d$x)=\sum\nolimits_{i\geq 1}e^{-x}\nu (-\log
s_{i}\in $d$x).$

Introduce then $T$, the first time at which $U$ and $U^{\prime }$ do not
belong to the same fragment and call $m(T)$ (resp. $m^{\prime }(T)$) the
length of the fragment containing $U$ (resp. $U^{\prime }$) at that time.
Since $\nu $ does not lose mass during sudden dislocations, the masses $m(T)$%
, $m^{\prime }(T)$ are a.s. strictly positive. Let then, for $m>0$, $\tau
(m\cdot )$ denote the function $t\in \left[ 0,\infty \right) \mapsto \tau
(mt)$. Using the fragmentation property, one sees that $D_{\tau }=T+%
\widetilde{D}_{\tau (m(T)\cdot )}$ and $D_{\tau }^{\prime }=T+\widetilde{D}%
_{\tau (m^{\prime }(T)\cdot )},$ where, conditionally on $m(T)$ and $%
m^{\prime }(T)$, $\widetilde{D}_{\tau (m(T)\cdot )}$ and $\widetilde{D}%
_{\tau (m^{\prime }(T)\cdot )}$\ are independent and distributed as $D_{\tau
(m(T)\cdot )}$ and $D_{\tau (m^{\prime }(T)\cdot )}$ respectively.
Therefore, $P(D_{\tau }=D_{\tau }^{\prime }<\infty )=P(\widetilde{D}_{\tau
(m(T)\cdot )}=\widetilde{D}_{\tau (m^{\prime }(T)\cdot )}<\infty )$ is equal
to $0$ as soon as the point $\infty $ is the only possible atom of $D_{\tau
(m\cdot )}$, $\forall m>0$. The proof ends with the following lemma. We
recall that $\sigma $ has no drift component.

\begin{lemma}
Let $f$ be a locally integrable and strictly positive function on $\left[
0,\infty \right) $. Suppose moreover that $f$ is monotone near $\infty $.
Then the integral $\int_{0}^{\infty }f(\sigma (r))\mathrm{d}r$ is either
a.s. finite or a.s. infinite and when it is a.s. finite, its distribution is
atomless.
\end{lemma}

\noindent \textbf{Proof. }The first assertion is a consequence of the
Hewitt-Savage $0$-$1$ law and is shown, e.g., in the proof of Prop.10 of %
\cite{Haas1}. In the following we suppose that the integral $%
\int_{0}^{\infty }f(\sigma (r))\mathrm{d}r$ is a.s. finite. In particular, $%
f $ is non-increasing near $\infty $ and converges to $0$.

(i). The proof is easy when $\nu $ is finite. Indeed, let then $T_{1}$ be
the first jump time of $\sigma .$ It is well-known that $T_{1}$ and $\sigma
(r+T_{1})$ are independent and that $T_{1}$ has an exponential distribution.
Therefore, splitting the integral at $T_{1}$, we see that $\int_{0}^{\infty
}f(\sigma (r))\mathrm{d}r$ can be written as the sum of two independent r.v.:%
\begin{equation*}
\int_{0}^{\infty }f(\sigma (r))\mathrm{d}r=f(0)T_{1}+\int_{0}^{\infty
}f(\sigma (r+T_{1}))\mathrm{d}r\text{,}
\end{equation*}%
the first one, $f(0)T_{1}$, being absolutely continuous. It is easy that $%
\int_{0}^{\infty }f(\sigma (r))\mathrm{d}r$ is then also absolutely
continuous, hence atomless.

(ii). From now on, we suppose that $\nu $ is infinite. Introduce then for
all $t\geq 0$ the stopping times
\begin{equation*}
\theta (t):=\inf \left\{ u:\int_{0}^{u}f(\sigma (r))\mathrm{d}r>t\right\} ,
\end{equation*}%
with the convention $\inf \{\emptyset \}=\infty $. According to the strong
Markov property, conditional on $\theta (t)<\infty ,$
\begin{equation*}
\int_{0}^{\infty }f(\sigma (r))\mathrm{d}r=t+\int_{0}^{\infty }f(\sigma
(\theta (t))+\sigma ^{(t)}(r))\mathrm{d}r
\end{equation*}%
where $\sigma ^{(t)}(r):=\sigma (r+\theta (t))-\sigma (\theta (t))$, $r\geq
0 $, is a subordinator distributed as $\sigma $ and independent of $\sigma
(\theta (t))$.

Now fix some $t>0$ and to begin with, suppose that $f$ is strictly
decreasing on $\left[ 0,\infty \right) $. The function $x\in \left( 0,\infty
\right) \mapsto \int_{0}^{\infty }f(x+\sigma ^{(t)}(r))\mathrm{d}r$ is then
strictly decreasing, hence injective.\ Consequently, there is at most one
point, say $X_{t}$, such that $\int_{0}^{\infty }f(X_{t}+\sigma ^{(t)}(r))%
\mathrm{d}r=t.$ If that point does not exist, $X_{t}:=\infty $. Then,%
\begin{eqnarray*}
P\left( \int_{0}^{\infty }f(\sigma (r))\mathrm{d}r=2t\right) &=&P\left(
\int_{0}^{\infty }f(\sigma (\theta (t))+\sigma ^{(t)}(r))\mathrm{d}%
r=t,\theta (t)<\infty \right) \\
&=&P(\sigma (\theta (t))=X_{t},\theta (t)<\infty )
\end{eqnarray*}%
with $X_{t}$ independent of $\sigma (\theta (t)).$ This latter probability
is then $0$, because for all $0<a<\infty $, $P(\sigma (\theta (t))=a)\leq
P(\exists s:\sigma (s)=a)$ and, by Kesten's Theorem (see e.g. Prop. 1.9 in %
\cite{Bertoin3}), since $\sigma $ has $0$ drift and $\nu $ is infinite, $%
P(\exists s:\sigma (s)=a)=0$. Hence the conclusion when $f$ is strictly
decreasing on $\left[ 0,\infty \right) $.

Suppose next that $f$ is only non-increasing on $\left[ 0,\infty \right) $
and that $P\left( \int_{0}^{\infty }f(\sigma (r))\mathrm{d}r=2t\right) >0$
for some $t>0.$ Still because $\sigma (\theta (t))$ has no atom (except $%
\infty $), this implies that the probability that the function $x\mapsto
\int_{0}^{\infty }f(x+\sigma ^{(t)}(r))\mathrm{d}r$ is constant on some
non-void open interval is strictly positive, which, because of the
monoticity of $f$, implies in turn that
\begin{equation*}
\exists \text{ }0<x<x^{\prime }:P\left( \int_{0}^{\infty }f(x+\sigma (r))%
\mathrm{d}r=\int_{0}^{\infty }f(x^{\prime }+\sigma (r))\mathrm{d}r=t\right)
>0.
\end{equation*}%
But this is not possible. Indeed, consider some sequence $t_{n}<t$, $%
t_{n}\rightarrow t$, such that $P\left( \int_{0}^{\infty }f(x^{\prime
}+\sigma (r))\mathrm{d}r=t_{n}\right) =0$, $\forall n.$ Then for all $n$,%
\begin{eqnarray*}
&&P\left( \int_{0}^{\infty }f(x+\sigma (r))\mathrm{d}r=\int_{0}^{\infty
}f(x^{\prime }+\sigma (r))\mathrm{d}r=t\right) \\
&=&P\left( \int_{0}^{\infty }f(y+\sigma (\theta (t-t_{n}))+\sigma
^{(t-t_{n})}(r))\mathrm{d}r=t_{n},\text{ }\forall y\in \left[ x,x^{\prime }%
\right] \right) .
\end{eqnarray*}%
The $t_{n}$'s have been chosen such that $\int_{0}^{\infty }f(x^{\prime
}+\sigma ^{(t-t_{n})}(r))\mathrm{d}r\neq t_{n}$ a.s. Therefore the latter
probability is necessarily smaller than $P\left( x^{\prime }<x+\sigma
(\theta (t-t_{n})\right) $, which tends to $0$ as $t_{n}\rightarrow t.$
Hence $P\left( \int_{0}^{\infty }f(\sigma (r))\mathrm{d}r=2t\right) =0$.

At last, when $f$ is non-increasing (only) in a neighborhood of $\infty ,$
say on $[b,\infty )$, we can turn down to the previous case as follows: let $%
T_{b}:=\inf \{t:\sigma (t)>b\}$ and write%
\begin{equation}
\int_{0}^{\infty }f(\sigma (r))\mathrm{d}r=\int_{0}^{T_{b}}f(\sigma (r))%
\mathrm{d}r+\int_{0}^{\infty }f(\sigma (T_{b})+\widetilde{\sigma }(r))%
\mathrm{d}r  \label{26}
\end{equation}%
where $\widetilde{\sigma }$ is independent of $(\sigma (t),t\leq T_{b})$ and
distributed as $\sigma .$ Conditional on $(\sigma (t),t\leq T_{b}),$ we know
that $\int_{0}^{\infty }f(\sigma (T_{b})+\widetilde{\sigma }(r))\mathrm{d}r$
is atomless since $f(\sigma (T_{b})+\cdot )$ is non-increasing on $\left[
0,\infty \right) $. Therefore, using (\ref{26}) and still conditioning on $%
(\sigma (t),t\leq T_{b}),$ we see that $\int_{0}^{\infty }f(\sigma (r))%
\mathrm{d}r$ is also atomless. \ \rule{0.5em}{0.5em}

\subsection{\textbf{Fragmentations with erosion}}

Until now, we have considered pure-jump fragmentation processes. However it
is well-known that a fragmentation may have a continuous part, and more
precisely, that a general homogeneous fragmentation is characterized by its
dislocation measure $\nu $ and by an \textit{erosion coefficient} $c\geq 0$
that measures the melting of the particles.\ More precisely, any homogeneous
fragmentation $F^{\text{hom}}$ can be factorized as $F^{\text{hom}%
}(t)=e^{-ct}\overline{F}^{\text{hom}}(t)$, $t\geq 0$, for some $c\geq 0$ and
some pure-jump $\nu $-homogeneous fragmentation $\overline{F}^{\text{hom}}$.
Exactly as in Section 1.1.1, one can then construct from any $(c,\nu )$%
-homogeneous fragmentation, some $(\tau ,c,\nu )$ fragmentation and $(\tau
,c,\nu ,I)$ fragmentation with immigration.

We still work under the hypothesis (\ref{H}). Theorems \ref{Thprincipal} and %
\ref{Proposition_infinie} can then be modified as follows:

- all the results concerning the convergence of $%
(F_{2}^{(m)},F_{3}^{(m)},...)$ are still valid, provided that in Theorem \ref%
{Thprincipal} we replace the $(\tau ,\nu ,I)$ fragmentation with immigration
by some $(\tau ,c,\nu ,I)$ fragmentation with immigration

- under the assumptions of Theorem \ref{Thprincipal}, this convergence holds
jointly with that of $(m-F_{1}^{(m)})/m\tau (m)$ to the deterministic
process $(ct,t\geq 0)$. Under the assumptions of Theorem \ref%
{Proposition_infinie}, it holds jointly with that of $(m-F_{1}^{(m)}((%
\varphi _{\nu }/\tau )(m)\cdot ))/m\varphi _{\nu }(m)\ $to $(ct,t\geq 0)$.

The main difference in the proofs is that the subordinator $\xi $ introduced
in $(\ref{31})$ is here replaced by the subordinator $\xi _{c}$, $\xi
_{c}(t):=\xi (t)+ct$, $t\geq 0$.


\begin{thebibliography}{99}
\bibitem{AldousCRT1} \textsc{D. Aldous}, The continuum random tree I,\emph{\
Annals of Probab}., \textbf{19} (1) (1991), pp. 1-28.

\bibitem{Aldous CRT3} \textsc{D. Aldous}, The continuum random tree III,%
\emph{\ Annals of Probab}., \textbf{21 }(1993), pp. 248-289.

\bibitem{AldousPitman} \textsc{D.~Aldous, J. Pitman}, The standard additive
coalescent, \emph{Ann. Probab.,} \textbf{26} (4) (1998), pp. 1703-1726.

\bibitem{berest02} \textsc{J.~Berestycki}, Ranked fragmentations,\emph{\
ESAIM Probab. Statist.}, \textbf{6} (2002), pp.~157--175.

\bibitem{BerestyckiThese} \textsc{J.~Berestycki}, \emph{Fragmentations et
coalescences homog\`{e}nes}, Th\`{e}se de doctorat de l'universit\'{e} Paris
6. Available via \texttt{http://tel.ccsd.cnrs.fr/}

\bibitem{Bertoin3} \textsc{J.~Bertoin}, Subordinators : Examples and
Applications. In P.\ Bernard (editor): \emph{Lectures on Probability Theory
and Statistics,}\textit{\ }Ecole d'\'{e}t\'{e} de probabilit\'{e}s de
St-Flour XXVII, pp. 1-91. Lect. Notes in Maths \textbf{1717}, Springer,
Berlin 1999.

\bibitem{berthfrag01} \textsc{J.~Bertoin}, Homogeneous fragmentation
processes, \emph{Probab. Theory Relat. Fields}, \textbf{121} (3) (2001),
pp.~301--318.

\bibitem{bertsfrag02} \textsc{J.~Bertoin}, Self-similar fragmentations,
\emph{Ann. Inst. Henri Poincar\'{e} Probab. Stat.}, \textbf{38} (2002),
pp.~319--340.

\bibitem{bertafrag02} \textsc{J.~Bertoin}, The asymptotic behavior of
fragmentation processes, \emph{J. Eur. Math. Soc.}, \textbf{5} (4) (2003),
pp. 395-416.

\bibitem{bgt} \textsc{N.H. Bingham, C.M. Goldie, J.L. Teugels}, \emph{%
Regular variation}, vol.~27 of Encyclopedia of Mathematics and its
Applications, Cambridge University Press, Cambridge 1989.

\bibitem{DuquesneImmigration} \textsc{T. Duquesne}, Continuum random trees
and branching processes with immigration, in preparation.

\bibitem{DuquesneLegall} \textsc{T. Duquesne, J.F. Le Gall}, \emph{Random
Trees, L\'{e}vy Processes and Spatial Branching Processes}, Ast\'{e}risque
281, Soci\'{e}t\'{e} Math\'{e}matique de France, 2002.

\bibitem{EthridgeWilliams} \textsc{A.M. Etheridge, D.E. Williams}, A
decomposition of the $1+\beta $ superprocess conditioned on survival, \emph{%
Proc. R.\ Soc. Edin.}, A \textbf{133} (2003), pp. 829-847.

\bibitem{Ethier Kurtz} \textsc{S.N. Ethier, T.G. Kurtz}, \emph{Markov
Processes, Characterization and Convergence}, Wiley and Sons, New-York 1986.

\bibitem{Evans} \textsc{S.N. Evans}, Two representations of a conditioned
superprocess, \emph{Proc. R.\ Soc. Edin.}, A \textbf{123} (1993), pp.
959-971.

\bibitem{Gromov} \textsc{M. Gromov}, \emph{Metric Structures for Riemannian
and Non-Riemannian Spaces. Progress in Mathematics. }Birkh\"{a}user, Boston
1999.

\bibitem{Haas1} \textsc{B. Haas}, Loss of mass in deterministic and random
fragmentations, \emph{Stoch. Proc. App.},\textit{\ }\textbf{106 }(2) (2003),
pp. 245-277.

\bibitem{HaasImmig04} \textsc{B. Haas}, Equilibrium for fragmentations with
immigration, \emph{Ann. Appl. Probab.}, \textbf{15 }(3) (2005) pp. 1958-1996.

\bibitem{HaasThese} \textsc{B. Haas}, \emph{Fragmentations et perte de masse}%
, Th\`{e}se de doctorat de l'universit\'{e} Paris 6. Available via \texttt{%
http://tel.ccsd.cnrs.fr/}

\bibitem{HaasMiermont} \textsc{B. Haas, G. Miermont}, The genealogy of
self-similar fragmentations with a negative index as a continuum random
tree, \emph{Elect. J. Probab.}, \textbf{9} (2004), pp. 57-97.

\bibitem{JacodShiryaev} \textsc{J. Jacod, A.N. Shiryaev}, \emph{Limit
Theorems for stochastic Processes}, second edition. Springer-Verlag, Berlin
2003.

\bibitem{Kallenberg} \textsc{O. Kallenberg}, \emph{Random Measures},
Akademie-Verlag, Berlin 1975.

\bibitem{kingman93} \textsc{J.F.C.~Kingman}, \emph{Poisson processes},
vol.~3 of Oxford Studies in Probability, The Clarendon Press Oxford
University Press, New York 1993.

\bibitem{Lambert} \textsc{A.~Lambert}, The genealogy of continuous-state
branching processes with immigration, \emph{Probab. Theory Relat. Fields,}
\textbf{122} (1) (2002), pp. 42--70.

\bibitem{Lyons Pemantle Peres} \textsc{R. Lyons, R. Pemantle, Y. Peres},
Conceptual proofs of \textit{L}Log\textit{L} criteria for mean behavior of
branching processes, \emph{Annals of Probab}., \textbf{23 }(3) (1995), pp.
1125-1138.

\bibitem{mierfmoins} \textsc{G.~Miermont}, Self-similar fragmentations
derived from the stable tree I: splitting at heights, \emph{Probab. Theory
Relat. Fields}, \textbf{127} (3) (2003), pp. 423-454.

\bibitem{mierfplus} \textsc{G.~Miermont}, Self-similar fragmentations
derived from the stable tree II: splitting at nodes, \emph{Probab. Theory
Relat. Fields}, \textbf{131} (3) (2005), pp. 341-375.

\bibitem{MiermontSchweinsberg} \textsc{G.~Miermont, J. Schweinsberg},
Self-similar fragmentations and stable subordinators, \emph{S\'{e}minaire de
Probabilit\'{e}s XXXVII, Lectures Notes in Math., 1832},\emph{\ }pp.
333-359, Springer, Berlin 2003.

\bibitem{Neveu} \textsc{J. Neveu}, Processus ponctuels. Ecole d'\'{e}t\'{e}
de probabilit\'{e}s de St-Flour VI, pp. 249-445. Lect. Notes in Maths
\textbf{598}, Springer, Berlin 1977.

\bibitem{Pitman} \textsc{J.~Pitman}, Combinatorial stochastic processes.
Lecture notes for St. Flour course, July 2002. To appear. Available via
\texttt{www.stat.berkeleley.edu.}
\end{thebibliography}
\end{document}